\documentclass{SCAE}
\numberwithin{equation}{section}

\def\ep{\textsf{E}} 
\def\Sbep{\widehat{\mathbb E}} 
\def\cSbep{\widehat{\mathcal E}} 
\def\extSbep{\mathbb E} 
\def\Capc{\mathbb V} 
\def\cCapc{\mathcal V} 
\def\outCapc{\mathbb V^{\ast}}
\def\outcCapc{\mathcal V^{\ast}}

\begin{document}

\Year{2016} %
\Month{Month}
\Vol{??} %
\No{?} %
\BeginPage{1} %
\EndPage{XX} %
\AuthorMark{Zhang L X}
\ReceivedDay{Month Date, Year}
\AcceptedDay{Month Date, Year}
\PublishedOnlineDay{; published online Month Date, Year}
\DOI{10.1007/s11425-000-0000-0} 

\title[Inequalities  under sub-linear expectations]{Exponential  inequalities  under the sub-linear expectations with applications to laws of the iterated logarithm}{}


\author[1]{ZHANG Li-Xin}{}

\address[{\rm1}]{School of Mathematical Sciences, Zhejiang University, Hangzhou {\rm 310027}, P.R. China;}

\Emails{ stazlx@zju.edu.cn}\maketitle


 {\begin{center}
\parbox{14.5cm}{\begin{abstract}
Kolmogorov's exponential  inequalities  are basic tools for studying  the strong limit theorems such as the classical laws of the iterated logarithm
for both independent and dependent random variables.
This paper establishes the Kolmogorov type exponential
inequalities of the  partial sums of independent random variables  as well as negatively dependent random variables under the sub-linear expectations.
As applications of the exponential inequalities,   the laws of the iterated logarithm in the sense of non-additive capacities are proved
for independent or negatively dependent identically distributed random variables  with finite second order moments.
For deriving  a lower bound of an exponential inequality,
 a central limit theorem is also proved under the sub-linear expectation for  random variables  with only finite  variances.
\vspace{-3mm}
\end{abstract}}\end{center}}



\renewcommand{\baselinestretch}{1.2}
\begin{center} \renewcommand{\arraystretch}{1.5}
{\begin{tabular}{lp{0.8\textwidth}} \hline \scriptsize
{\bf Citation:}\!\!\!\!&\scriptsize Zhang L X. \makeatletter\@titlehead.
Sci China Math, 2016, 59,
 doi:~\@DOI\makeatother\vspace{1mm}
\\
\hline
\end{tabular}}\end{center}

\baselineskip 11pt\parindent=10.8pt  \wuhao

\section{ Introduction and notations.}\label{sect1}
\setcounter{equation}{0}

Non-additive probabilities and non-additive expectations  are useful tools for studying  uncertainties in statistics, measures of risk,
superhedging in finance and non-linear stochastic calculus, cf. \cite{DM06}\cite{Gilboa87}\cite{Marinacci99}\cite{Peng97}-\cite{Peng08a} etc.
This paper considers the  general sub-linear expectations and related non-additive probabilities generated by them.
The general framework of the sub-linear expectation is introduced by Peng \cite{Peng06}-\cite{Peng08b} in a general function space by relaxing   the linear property of the classical  expectation
to  the sub-additivity
and positive homogeneity (cf. Definition~\ref{def1.1} below).
The sub-linear expectation  provides
a  very flexible framework to model non-additive probability problems  and  produces many interesting properties different
 from those of the linear expectations.
  For example, one constant is not enough to characterize
  the mean or variance of a random variable  in a sub-linear expectation space, the limit in the law of large numbers is no longer a contact,
  and, comparing to the  classical one-dimensional  normal  distribution which is characterized by the Stein equation, an
ordinary differential equation (ODE),
a normal distribution under the sub-linear expectation is characterized by a time-space parabolic partial  different equation (PDE). Recently, Hu and Li \cite{HuLi14} showed that the characteristic function cannot determine the distribution of
random variables on the sub-linear expectation space.
  Roughly speaking, a sub-linear expectation is related to a group of unknown linear expectations
and the distribution under a sub-linear expectation is related to a group of probabilities (cf. Lemma 2.4 of Peng \cite{Peng08b}).
For more properties of the sub-linear expectations, one can refer to Peng \cite{Peng08b},
where the notion of independent and identically distributed random variables under  the sub-linear expectations is introduced  and
the weak convergence such as  central limit theorems and weak laws of large numbers are studied.

The motivation  of this paper is to study the laws of the iterated logarithm under reasonable conditions. A law of the iterated logarithm for independent and identically distributed random variables in the  sub-linear expectation space was established recently by Chen and Hu \cite{ChenHu14}. However, a very strict condition that the random variables are bounded is assumed.
Basically,   the classical law of the iterated logarithm is established through
the Kolmogorov type exponential inequalities for both independent and negatively dependent random variables (cf. Petrov \cite{Petrov95}, Shao and Su \cite{ShaoSu99}).
The main purpose of this paper is to establish the Kolmogorov type exponential inequalities  for   independent random variables as well as negatively dependent random variables
    in the general sub-linear expectation spaces.  By applying these inequalities,
    we prove that the laws of the iterated logarithm holds for independent random variables as well as negatively dependent   random variables
    under the condition that only the second order moments are finite.   It is shown that for a sequence $\{X_n;n\ge 1\}$ of
    independent and identically distributed random variables with finite variances, the law of the iterated logarithm holds if and only if
    the sub-linear means  are zeros and  the  Choquet integral    of $X_1^2/\log\log |X_1|$ is finite.  Also, for deriving   a lower bound of an exponential inequality    for
   independent and identically distributed random variables,   we   prove a central limit theorem under
   only the condition that the second order moments are finite,
   which improves the central limit theorem of Peng \cite{Peng08b} (cf. Remark~\ref{remarkCLT} below).
    Because the sub-linear expectation is not additive, many powerful tools for linear expectations and probabilities
    such as the martingale method, the stopping time, the  symmetrization method are not valid, so that the study of the limit theorems becomes much more technical even after the  exponential inequalities are established.
    In the next section, we give some notations under the sub-linear expectations including independence and negative dependence.
   In    Section~\ref{sectMain}, we give the main results. The proof is given in the last section.

\section{Basic Settings}\label{sectBasic}
\setcounter{equation}{0}

We use the framework and notations of Peng \cite{Peng08b}. Let  $(\Omega,\mathcal F)$
 be a given measurable space  and let $\mathscr{H}$ be a linear space of real functions
defined on $(\Omega,\mathcal F)$ such that if $X_1,\ldots, X_n \in \mathscr{H}$  then $\varphi(X_1,\ldots,X_n)\in \mathscr{H}$ for each
$\varphi\in C_{l,Lip}(\mathbb R_n)$,  where $C_{l,Lip}(\mathbb R_n)$ denotes the linear space of (local Lipschitz)
functions $\varphi$ satisfying
\begin{eqnarray*} & |\varphi(\bm x) - \varphi(\bm y)| \le  C(1 + |\bm x|^m + |\bm y|^m)|\bm x- \bm y|, \;\; \forall \bm x, \bm y \in \mathbb R_n,&\\
& \text {for some }  C > 0, m \in \mathbb  N \text{ depending on } \varphi. &
\end{eqnarray*}
$\mathscr{H}$ is considered as a space of ``random variables''. In this case we denote $X\in \mathscr{H}$.

\begin{definition}\label{def1.1} A  sub-linear expectation $\Sbep$ on $\mathscr{H}$  is a function $\Sbep: \mathscr{H}\to \overline{\mathbb R}$ satisfying the following properties: for all $X, Y \in \mathscr H$, we have
\begin{description}
  \item[\rm (a)]  Monotonicity: If $X \ge  Y$ then $\Sbep [X]\ge \Sbep [Y]$;
\item[\rm (b)] Constant preserving: $\Sbep [c] = c$;
\item[\rm (c)] Sub-additivity: $\Sbep[X+Y]\le \Sbep [X] +\Sbep [Y ]$ whenever $\Sbep [X] +\Sbep [Y ]$ is not of the form $+\infty-\infty$ or $-\infty+\infty$;
\item[\rm (d)] Positive homogeneity: $\Sbep [\lambda X] = \lambda \Sbep  [X]$, $\lambda\ge 0$.
 \end{description}
 Here $\overline{\mathbb R}=[-\infty, \infty]$. The triple $(\Omega, \mathscr{H}, \Sbep)$ is called a sub-linear expectation space. Give a sub-linear expectation $\Sbep $, let us denote the conjugate expectation $\cSbep$of $\Sbep$ by
$$ \cSbep[X]:=-\Sbep[-X], \;\; \forall X\in \mathscr{H}. $$
\end{definition}

From the definition, it is easily shown that    $\cSbep[X]\le \Sbep[X]$, $\Sbep[X+c]= \Sbep[X]+c$ and $\Sbep[X-Y]\ge \Sbep[X]-\Sbep[Y]$ for all
$X, Y\in \mathscr{H}$ with $\Sbep[Y]$ being finite. Further, if $\Sbep[|X|]$ is finite, then $\cSbep[X]$ and $\Sbep[X]$ are both finite.  Denote
$ \mathscr{L}=\{X\in \mathscr{H}: \Sbep[|X|]<\infty\}. $

\begin{definition} ({\em See \cite{Peng06}\cite{Peng08b}})

\begin{description}
  \item[ \rm (i)] ({\em Identical distribution}) Let $\bm X_1$ and $\bm X_2$ be two $n$-dimensional random vectors defined
respectively in sub-linear expectation spaces $(\Omega_1, \mathscr{H}_1, \Sbep_1)$
  and $(\Omega_2, \mathscr{H}_2, \Sbep_2)$. They are called identically distributed, denoted by $\bm X_1\overset{d}= \bm X_2$  if
$$ \Sbep_1[\varphi(\bm X_1)]=\Sbep_2[\varphi(\bm X_2)], \;\; \forall \varphi\in C_{l,Lip}(\mathbb R_n), $$
whenever the sub-expectations are finite. A sequence $\{X_n;n\ge 1\}$ of random variables is said to be identically distributed if $X_i\overset{d}= X_1$ for each $i\ge 1$.
\item[\rm (ii)] ({\em Independence})   In a sub-linear expectation space  $(\Omega, \mathscr{H}, \Sbep)$, a random vector $\bm Y =
(Y_1, \ldots, Y_n)$, $Y_i \in \mathscr{H}$ is said to be independent to another random vector $\bm X =
(X_1, \ldots, X_m)$ , $X_i \in \mathscr{H}$ under $\Sbep$  if for each test function $\varphi\in C_{l,Lip}(\mathbb R_m \times \mathbb R_n)$
we have
$ \Sbep [\varphi(\bm X, \bm Y )] = \Sbep \big[\Sbep[\varphi(\bm x, \bm Y )]\big|_{\bm x=\bm X}\big],$
whenever $\overline{\varphi}(\bm x):=\Sbep\left[|\varphi(\bm x, \bm Y )|\right]<\infty$ for all $\bm x$ and
 $\Sbep\left[|\overline{\varphi}(\bm X)|\right]<\infty$.
 \item[\rm (iii)] ({\em IID random variables}) A sequence of random variables $\{X_n; n\ge 1\}$
 is said to be independent and identically distributed (IID), if
 $X_i\overset{d}=X_1$ and $X_{i+1}$ is independent to $(X_{i+1},\ldots, X_n)$ for each $i\ge 1$.
 \end{description}
\end{definition}

From the definition of independence, it is easily seen that, if $Y$ is independent to $X$, and $X, Y\in \mathscr{L}$,  $X\ge 0, \Sbep [Y]\ge 0$, then
\begin{equation}\label{eq1.1} \Sbep[XY]=\Sbep[X]\Sbep[Y].
\end{equation}
Further,  if $Y$ is independent to $X$ and $0\le X, Y\in \mathscr{L}$, then
\begin{equation}\label{eq1.2} \Sbep[XY]=\Sbep[X]\Sbep[Y], \;\; \cSbep[XY]=\cSbep[X]\cSbep[Y].
\end{equation}

Motivated by the above properties (\ref{eq1.1}) and (\ref{eq1.2}),    we give the concept of negative dependence under the sub-linear expectation.

 \begin{definition} \begin{description}
 \item[\rm (i)]  ({\em Negative dependence})  In a sub-linear expectation space  $(\Omega, \mathscr{H}, \Sbep)$, a random vector $\bm Y =
(Y_1, \ldots, Y_n)$, $Y_i \in \mathscr{H}$ is said to be negatively dependent (ND) to another random vector $\bm X =
(X_1, \ldots, X_m)$, $X_i \in \mathscr{H}$ under $\Sbep$  if for each pair of   test functions $\varphi_1\in C_{l,Lip}(\mathbb R_m)$ and
$\varphi_2\in C_{l,Lip}(\mathbb R_n)$
we have
$ \Sbep [\varphi_1(\bm X)\varphi_2(\bm Y )] \le  \Sbep [\varphi_1(\bm X)]\Sbep[\varphi_2(\bm Y )]$
whenever either $\varphi_1, \varphi_2$ are coordinatewise nondecreasing or $\varphi_1, \varphi_2$ are coordinatewise non-increasing with $\varphi_1(\bm X)\ge 0$, $\Sbep[\varphi_2(\bm Y )]\ge 0$, $\Sbep [|\varphi_1(\bm X)\varphi_2(\bm Y )|]<\infty$,
$\Sbep [|\varphi_1(\bm X)|]<\infty$,
  $\Sbep [|\varphi_2(\bm Y )|]<\infty$.
  \item[\rm (ii)] ({\em ND random variables})   Let $\{X_n;n\ge 1\}$ be a sequence of random variables in the sub-linear expectation space
  $(\Omega, \mathscr{H}, \Sbep)$. $X_1,X_2,\ldots $ are said to   be negatively dependent if $X_{i+1}$ is negatively dependent
to $(X_1,\ldots, X_i)$ for each $i\ge 1$.
\end{description}
\end{definition}
It is obvious that, if $\{X_n;n\ge 1\}$ is a sequence of  independent random variables and $f_1(x),f_2(x),\ldots\in C_{l,Lip}(\mathbb R)$,
 then $\{f_n(X_n);n\ge 1\}$ is  also a sequence of independent random variables;  if $\{X_n;n\ge 1\}$ is  a sequence of negatively dependent random variables and $f_1(x),f_2(x),\ldots\in C_{l,Lip}(\mathbb R)$
 are non-decreasing (resp.  non-increasing) functions, then $\{f_n(X_n);n\ge 1\}$ is  also a sequence of negatively dependent random variables.

Next, we introduce the capacities corresponding to the sub-linear expectations.
Let $\mathcal G\subset\mathcal F$. A function $V:\mathcal G\to [0,1]$ is called a capacity if
$$ V(\emptyset)=0, \;V(\Omega)=1 \; \text{ and } V(A)\le V(B)\;\; \forall\; A\subset B, \; A,B\in \mathcal G. $$
It is called to be sub-additive if $V(A\bigcup B)\le V(A)+V(B)$ for all $A,B\in \mathcal G$  with $A\bigcup B\in \mathcal G$.

 Let $(\Omega, \mathscr{H}, \Sbep)$ be a sub-linear space,
and  $\cSbep $  be  the conjugate expectation of $\Sbep$.  We denote a pair $(\Capc,\cCapc)$ of capacities by
$$ \Capc(A):=\inf\{\Sbep[\xi]: I_A\le \xi, \xi\in\mathscr{H}\}, \;\; \cCapc(A):= 1-\Capc(A^c),\;\; \forall A\in \mathcal F, $$
where $A^c$  is the complement set of $A$.
Then
\begin{equation}\label{eq1.3} \begin{matrix}
&\Capc(A):=\Sbep[I_A], \;\; \cCapc(A):= \cSbep[I_A],\;\; \text{ if } I_A\in \mathscr H\\
&\Sbep[f]\le \Capc(A)\le \Sbep[g], \;\;\cSbep[f]\le \cCapc(A) \le \cSbep[g],\;\;
\text{ if } f\le I_A\le g, f,g \in \mathscr{H}.
\end{matrix}
\end{equation}
It is obvious that $\Capc$ is sub-additive. But $\cCapc$ and $\cSbep$ are not. However, we have
\begin{equation}\label{eq1.4}
  \cCapc(A\bigcup B)\le \cCapc(A)+\Capc(B) \;\;\text{ and }\;\; \cSbep[X+Y]\le \cSbep[X]+\Sbep[Y]
\end{equation}
due to the fact that $\Capc(A^c\bigcap B^c)=\Capc(A^c\backslash B)\ge \Capc(A^c)-\Capc(B)$ and $\Sbep[-X-Y]\ge \Sbep[-X]-\Sbep[Y]$.

Also, we define the  Choquet integrals/expecations $(C_{\Capc},C_{\cCapc})$  by
$$ C_V[X]=\int_0^{\infty} V(X\ge t)dt +\int_{-\infty}^0\left[V(X\ge t)-1\right]dt $$
with $V$ being replaced by $\Capc$ and $\cCapc$ respectively.
It can be verified that (cf., Lemma \ref{lem2} (iii)), if $\lim_{c\to \infty}\Sbep[(|X|-c)^+]=0$, then
$ \Sbep[|X|]\le C_{\Capc}(|X|).$

\section{Main results}\label{sectMain}
\setcounter{equation}{0}
In this section, we give the mains results.
We first  give the upper bounds of the  exponential inequalities for independent random variables as well as   negatively dependent random variables,
then   a lower  bound of an exponential inequality for independent and identically distributed random variables.
 For deriving this lower bound, we give a new central limit theorem. At last, we give the laws of the iterated logarithm.

 \subsection{Exponential  inequalities} \label{sectIneq}

Let $\{X_1,\ldots, X_n\}$ be a sequence  of   random variables
in $(\Omega, \mathscr{H}, \Sbep)$ . Set $S_n=\sum_{k=1}^n X_k$,  $B_n=\sum_{k=1}^n \Sbep[X_k^2]$ and $M_{n,p}=\sum_{k=1}^n\Sbep[|X_k|^p]$, $p\ge 2$.
The following is  our main result on the exponential inequalities.
For the exponential inequalities for classical negatively dependent random variables, one can refer to Su, Zhao and Wang \cite{SZW97}, Shao \cite{Shao00} etc.
The comparison method of Shao \cite{Shao00} is not valid under the sub-linear expectation $ \Sbep$  because of the non-additivity.
\begin{theorem}\label{th2}    Let $\{X_1,\ldots, X_n\}$ be a sequence  of negatively dependent random variables
in $(\Omega, \mathscr{H}, \Sbep)$ with $\Sbep[X_k]\le 0$. Then
\begin{description}
  \item[\rm (a)]
  For all $x,y>0$,
\begin{equation}\label{eqth2.1}
\Capc\left(S_n\ge x\right)\le \Capc\left(\max_{k\le n}X_k\ge y\right)
+ \exp\left\{-\frac{x^2}{2(xy+B_n)}\Big(1+\frac{2}{3}\ln \big(1+\frac{xy}{B_n}\big)\Big)\right\};
\end{equation}
\item[\rm (b)] For any $p\ge 2$, there exists a constant $C_p\ge 1$ such that
\begin{equation}\label{eqth2.2}
\Capc\left(S_n\ge x\right)\le C_p \delta^{-2p} \frac{M_{n,p}}{x^p}
+\exp\left\{-\frac{x^2}{2 B_n (1+\delta)}\right\}, \;\; \forall x>0 \; \text{ and } 0<\delta\le 1;
\end{equation}
\item[\rm (c)] We have
\begin{align}\label{eqth2.3}C_{\Capc}\left[(S_n^+)^p\right]
\le  & p^p C_{\Capc}\Big[\big(\max_{k\le n}X_k^+\big)^p\Big]+ C_p B_n^{p/2} \nonumber \\
\le &
p^p \sum_{k=1}^n C_{\Capc}\Big[(X_k^+)^p\Big]+ C_p B_n^{p/2}, \;\; \forall p\ge 2.
\end{align}
\end{description}
 \end{theorem}

The following corollary gives the estimates of $\cCapc\left(S_n\ge x\right)$.

\begin{corollary}\label{cor2}    Let $\{X_1,\ldots, X_n\}$ be a sequence  of independent random variables
in $(\Omega, \mathscr{H}, \Sbep)$ with $\cSbep[X_k]\le 0$. Then
\begin{description}
  \item[\rm (a)]
  For all $x,y>0$,
\begin{equation}\label{eqcor2.1}
\cCapc\left(S_n\ge x\right)\le \Capc\left(\max_{k\le n}X_k\ge y\right)
+ \exp\left\{-\frac{x^2}{2(xy+B_n)}\Big(1+\frac{2}{3}\ln \big(1+\frac{xy}{B_n}\big)\Big)\right\};
\end{equation}
\item[\rm (b)] For any $p\ge 2$, there exists a constant $C_p\ge 1$ such that
\begin{equation}\label{eqcor2.2}
\cCapc\left(S_n\ge x\right)\le C_p \delta^{-2p} \frac{M_{n,p}}{x^p}
+\exp\left\{-\frac{x^2}{2 B_n (1+\delta)}\right\}, \;\; \forall x>0 \; \text{ and } 0<\delta\le 1.
\end{equation}
\end{description}
 \end{corollary}
 By choosing $p=2$ and $\delta=1$ in (\ref{eqcor2.2}) and applying the inequality $xe^{-x}\le e^{-1}$ ($x\ge 0$) we obtain
 \begin{equation}\label{eqcor2.4}
\cCapc\left(S_n\ge x\right)\le C\frac{\sum_{k=1}^n \Sbep[X_k^2]}{x^2}, \; \forall x>0.
\end{equation}
\bigskip

 The next theorem  give a lower bound of an    exponential inequality  for independent  and identically distributed random variables.
\begin{theorem}\label{thLower} Suppose that $\{X_n; n\ge 1\}$ is a sequence of   independent  and identically distributed random variables with
   $\Sbep[X_1]=\Sbep[-X_1]=0$ and
   $\lim\limits_{c\to \infty}  \Sbep\left[(X_1^2-c)^+\right]=0$.
   Write $\overline{\sigma}^2=\Sbep[X_1^2]$ and $\underline{\sigma}^2=\cSbep[X_1^2]$. Let $\{y_n\}$ be a sequence of positive numbers such that
   $y_n\to \infty$, $y_n/\sqrt{n}\to 0$.  Then
   \begin{description}
     \item[\rm (a) ]
         for any  $|b|<\underline{\sigma}$, $\epsilon>0$ and $\delta>0$
   with $(b/\underline{\sigma})^2+\delta<1$, there exists $n_0$ such that
\begin{equation}\label{eqLowerB} \cCapc\left(\left|\frac{S_n}{y_n\sqrt{n} }-b\right|\le \epsilon\right)\ge \exp\left\{-\left(\Big(\frac{|b|}{\underline{\sigma}}\Big)^2+\delta\right)\frac{y_n^2}{2} \right\},
   \;\; \forall n\ge n_0;
   \end{equation}
   \item[\rm (b) ]
         for any  $|b|<\overline{\sigma}$, $\epsilon>0$ and $\delta>0$
   with $(b/\overline{\sigma})^2+\delta<1$, there exists $n_0$ such that
\begin{equation}\label{eqLowerB.2} \Capc\left(\left|\frac{S_n}{y_n\sqrt{n} }-b\right|\le \epsilon\right)\ge \exp\left\{-\left(\Big(\frac{|b|}{\overline{\sigma}}\Big)^2+\delta\right)\frac{y_n^2}{2} \right\},
   \;\; \forall n\ge n_0.
   \end{equation}
    \end{description}
\end{theorem}
\begin{remark} (\ref{eqth2.3}) is the  Rosenthal type inequality under the Choquet expectation. (\ref{eqth2.2}) and (\ref{eqcor2.2}) are the Fuk and Nagaev\cite{FukNagaev71} type inequalities.
 (\ref{eqth2.1}) and (\ref{eqcor2.1}) are  the upper bounds of the Kolmogorov type exponential  inequalities. (\ref{eqLowerB}) and (\ref{eqLowerB.2})  are the lower bounds (cf. Lemmas~7.1 and 7.2 of Petrov \cite{Petrov95}).
However, the lower bounds can not be established in the same as that of Lemma 7.2 of Petrov \cite{Petrov95} because the sub-linear expectation is not additive over unjoint events.
\end{remark}

\subsection{ A central limit theorem}\label{sectCLT}

For proving  Theorem \ref{thLower}, we need the following central limit theorem for independent and identically distributed random variables
with only finite variances, which improves the central limit theorem of Peng \cite{Peng08b} and is of independent interest.
 \begin{theorem}\label{thCLT} (CLT) Suppose that $\{X_n; n\ge 1\}$ is a sequence of   independent  and identically distributed random variables with
   $\Sbep[X_1]=\Sbep[-X_1]=0$ and
   $\lim\limits_{c\to \infty}  \Sbep\left[(X_1^2-c)^+\right]=0$.  Write $\overline{\sigma}^2=\Sbep[X_1^2]$ and $\underline{\sigma}^2=\cSbep[X_1^2]$. Then for any   continuous function $\varphi$ satisfying  $|\varphi(x)|\le C (1+x^2)$,
  \begin{equation} \label{eqCLT} \lim_{n\to \infty}\Sbep\left[\varphi\left(\frac{S_n}{\sqrt{n}}\right)\right]=\widetilde{\mathbb E}[\varphi(\xi)],
  \end{equation}
   where   $\xi\sim N(0,[\underline{\sigma}^2, \overline{\sigma}^2])$  under $\widetilde{\mathbb E}$.
    Further, if $p> 2$ and $\Sbep[|X_1|^p]<\infty$,
     then (\ref{eqCLT}) holds for any   continuous function $\varphi$ satisfying  $|\varphi(x)|\le C (1+|x|^p)$ .
 \end{theorem}
Here, a random variable $\xi$ in a sub-linear expectation space $(\widetilde{\Omega}, \widetilde{\mathscr H}, \widetilde{\mathbb E})$   is call a normal $N\big(0, [\underline{\sigma}^2, \overline{\sigma}^2]\big)$ distributed random variable  (write $X\sim N\big(0, [\underline{\sigma}^2, \overline{\sigma}^2]\big)$  under $\widetilde{\mathbb E}$), if for any bounded Lipschitz function $\varphi$, the function $u(x,t)=\widetilde{\mathbb E}\left[\varphi\left(x+\sqrt{t} \xi\right)\right]$ ($x\in \mathbb R, t\ge 0$) is the unique viscosity solution of  the following heat equation:
      $$ \partial_t u -G\left( \partial_{xx}^2 u\right) =0, \;\; u(0,x)=\varphi(x), $$
where $G(\alpha)=\frac{1}{2}(\overline{\sigma}^2 \alpha^+ - \underline{\sigma}^2 \alpha^-)$.

\begin{remark}\label{remarkCLT} Peng \cite{Peng08b} pointed  that (\ref{eqCLT}) holds for all continuous function $\varphi$ satisfying a polynomial growth condition:
 $|\varphi(x)|\le C (1+|x|^k)$ for some $k$ (cf. his Theorem 5.1).
 However, in his  proof   for a bounded and  Lipschitz continuous function $\varphi$,
 the $(2+\alpha)$-th moment $\Sbep[|X_i|^{2+\alpha}]$ needs to be assumed bounded (cf. the proof of his Lemma 5.4).
Also, when a continuous function is extended to a continuous function $\varphi$ satisfying  $|\varphi(x)|\le C (1+|x|^{p-1})$, the following condition
 is needed (cf. his Lemma~5.5, where $Y_n=S_n/\sqrt{n}$):
 $$ \sup_n\Sbep\left[\Big|\frac{S_n}{\sqrt{n}}\Big|^p\right]<\infty, $$
 which is not verified in Peng \cite{Peng08b}. Such moment inequalities are not obvious under the sub-linear expectations (cf. Zhang \cite{Zh15b}).
   Now, note  Theorem \ref{thCLT} and
   $$\Sbep[|X_1|^p] =\Sbep[|X_n|^p]\le 2^{p-1}\big(\Sbep[|S_n|^p]+\Sbep[|S_{n-1}|^{p}]\big). $$
    A sufficient and necessary condition for (\ref{eqCLT}) to hold for any continuous function satisfying a polynomial growth condition is that
$\Sbep[|X_1|^p]<\infty$ for all $p>0$. It is also important to note that $|x|^k$ is a continuous function  satisfying a polynomial growth condition, but
$e^x$ and $e^{x^2}$ are not.
\end{remark}

\begin{remark} When we prepare this paper, we establish a functional central limit by applying Theorem \ref{thCLT}. One can refer to Zhang \cite{Zh15a}, where  Chung's law of the iterated logarithm is also established.
\end{remark}
  \subsection{The law of the iterated logarithm}\label{sectLIL}

Before we give the laws of the iterated logarithm,  we need some more notations about the sub-linear expectations and capacities.

\begin{definition}\label{def3.1}
\begin{description}
\item{\rm (I)} A sub-linear expectation $\Sbep: \mathscr{H}\to \mathbb R$ is called to be  countably sub-additive if it satisfies
\begin{description}
  \item[\rm (e)] {\em Countable sub-additivity}: $\Sbep[X]\le \sum_{n=1}^{\infty} \Sbep [X_n]$, whenever $X\le \sum_{n=1}^{\infty}X_n$,
  $X, X_n\in \mathscr{H}$ and
  $X\ge 0, X_n\ge 0$, $n=1,2,\ldots$;
 \end{description}
It is called to be continuous if it satisfies
\begin{description}
  \item[\rm (f) ]  {\em Continuity from below}: $\Sbep[X_n]\uparrow \Sbep[X]$ if $0\le X_n\uparrow X$, where $X_n, X\in \mathscr{H}$;
  \item[\rm (g) ] {\em Continuity from above}: $\Sbep[X_n]\downarrow \Sbep[X]$ if $0\le X_n\downarrow X$, where $X_n, X\in \mathscr{H}$.
\end{description}

\item{\rm (II)}  A function $V:\mathcal F\to [0,1]$ is called to be  countably sub-additive if
$$ V\Big(\bigcup_{n=1}^{\infty} A_n\Big)\le \sum_{n=1}^{\infty}V(A_n) \;\; \forall A_n\in \mathcal F. $$

\item{\rm (III)}  A capacity $V:\mathcal F\to [0,1]$ is called a continuous capacity if it satisfies
\begin{description}
  \item[\rm (III1) ] {\em Continuity from below}: $V(A_n)\uparrow V(A)$ if $A_n\uparrow A$, where $A_n, A\in \mathcal F$;
  \item[\rm (III2) ] {\em Continuity from above}: $V(A_n)\downarrow  V(A)$ if $A_n\downarrow A$, where $A_n, A\in \mathcal F$.
\end{description}
\end{description}
\end{definition}
\bigskip

It is obvious that a continuous sub-additive capacity $V$ (resp. a sub-linear expectation $\Sbep$) is countably sub-additive.
The ``the convergence part'' of the   Borel-Cantelli Lemma is still true for a countably sub-additive capacity.
\begin{lemma} ({\em Borel-Cantelli's Lemma}) Let $\{A_n, n\ge 1\}$ be a sequence of events in $\mathcal F$.
Suppose that $V$ is a countably sub-additive capacity.   If $\sum_{n=1}^{\infty}V\left (A_n\right)<\infty$, then
$$ V\left (A_n\;\; i.o.\right)=0, \;\; \text{ where } \{A_n\;\; i.o.\}=\bigcap_{n=1}^{\infty}\bigcup_{i=n}^{\infty}A_i. $$
\end{lemma}
{\bf Proof.} By the monotonicity and the countable sub-additivity, it follows that
\begin{align*}
0\le V\left (\bigcap_{n=1}^{\infty}\bigcup_{i=n}^{\infty}A_i\right)
\le V\left (\bigcup_{i=n}^{\infty}A_i\right)\le \sum_{i=n}^{\infty}V\left (A_i\right)\to 0\; \text{ as } n\to\infty. \;\Box
\end{align*}

\begin{remark}
It is important to note that the condition that ``$X$ is independent to $Y$ under $\Sbep$'' does not implies   that
``$X$ is independent to $Y$ under $\Capc$'' because the indicator functions $I\{X\in A\}$ and $I\{Y\in B\}$ are not in $C_{l,Lip}(\mathbb R)$,
and also, ``$X$ is independent to $Y$ under $\Capc$'' does not implies   that
``$X$ is independent to $Y$ under $\Sbep$'' because $\Sbep$ is not an integral with respect to $\Capc$.
So, we have not ``the divergence part'' of the Borel-Cantelli Lemma.
\end{remark}

 Because $\Capc$ may be not countably sub-additive in general, we define an  outer capacity $\outCapc$ by
$$ \outCapc(A)=\inf\Big\{\sum_{n=1}^{\infty}\Capc(A_n): A\subset \bigcup_{n=1}^{\infty}A_n\Big\},\;\; \outcCapc(A)=1-\outCapc(A^c),\;\;\; A\in\mathcal F.$$
  Then it can be shown that $\outCapc(A)$ is a countably sub-additive capacity with $\outCapc(A)\le \Capc(A)$ and the following properties:
  \begin{description}
    \item[\rm (a*)] If $\Capc$ is countably sub-additive, then  $\outCapc\equiv\Capc$.
    \item[\rm (b*)] If $I_A\le g$, $g\in \mathscr{H}$, then $\outCapc(A)\le \Sbep[g]$. Further, if $\Sbep$ is countably sub-additive, then
   \begin{equation}\label{eq3.2}
    \Sbep[f]\le \outCapc(A)\le \Capc(A)\le \Sbep[g], \;\; \forall f\le I_A\le g, f,g\in \mathscr{H}.
    \end{equation}
    \item[\rm (c*)]  $\outCapc$ is the largest countably sub-additive capacity satisfying
    the property that $\outCapc(A)\le \Sbep[g]$ whenever $I_A\le g\in \mathscr{H}$, i.e.,
    if $V$ is also a countably sub-additive capacity satisfying $V(A)\le \Sbep[g]$ whenever $I_A\le g\in \mathscr{H}$, then $V(A)\le \outCapc(A)$.
    \end{description}

In this subsection, we let $\{X_n;n\ge 1\}$ be a sequence of identically distributed  random variables in $(\Omega,\mathscr{H},\Sbep)$.
 Denote $\underline{\sigma}^2=\cSbep[X_1^2]$,
 $\overline{\sigma}^2=\Sbep[X_1^2]$,  $a_n=\sqrt{2n\log\log n}$, where $\log x=\ln(x\vee e)$.
The following is the law of the iterated logarithm for independent random variables and negatively dependent random variables.
For the  law of the iterated logarithm for classical  negatively dependent random variables,
 one can refer to Shao and Su \cite{ShaoSu99}, Zhang \cite{Zh01a} etc.
   \begin{theorem} \label{thLIL}
  \begin{description}
  \item[\rm (a)]   Suppose that $X_1, X_2, \ldots$ are  negatively dependent  with
   $\Sbep[X_1]=\Sbep[-X_1]=0$,
   $\lim\limits_{c\to \infty}  \Sbep\left[(X_1^2-c)^+\right]=0$ and
    \begin{equation}\label{eqthLIL.1} C_{\Capc}\left[\frac{X_1^2}{\log\log|X_1|}\right]<\infty.
    \end{equation}
       Then
\begin{equation}\label{eqthLIL.2}
\outCapc\left(\Big\{\liminf_{n\to \infty}\frac{S_n}{a_n}< -\overline{\sigma}\Big\}
\bigcup \Big\{\limsup_{n\to \infty}\frac{S_n}{a_n}> \overline{\sigma}\Big\}\right)=0.
\end{equation}
\item[\rm (b)] Suppose that $X_1,X_2,\ldots$ are  independent,   $\outCapc$ is continuous and  $\Sbep$ is countably sub-additive.
If
\begin{equation}\label{eqthLIL.3}
\outCapc\left(   \limsup_{n\to \infty}\frac{|S_n|}{a_n}=+\infty \right)<1,
\end{equation}
then $\Sbep[X_1]=\Sbep[-X_1]=0$, and  (\ref{eqthLIL.1}) holds.
\item[\rm (c)] Suppose that $X_1,X_2,\ldots$ are  independent,   $\Capc$ is continuous. Assume that (\ref{eqthLIL.3}) holds. Then we have  (\ref{eqthLIL.1}).
Further, if $\lim\limits_{c\to\infty}\Sbep[(|X_1|-c)^+]=0$, then  $\Sbep[X_1]=\Sbep[-X_1]=0 $.
\end{description}
\end{theorem}

\begin{remark} Theorem~\ref{thLIL} (a) can be regarded as the direct part which gives sufficient conditions for the law of the iterated logarithm
to hold,
and Theorem~\ref{thLIL} (b) (c) can be regarded the inverse part which gives the necessary conditions. According to this Theorem,
 we conjecture that the necessary and sufficient conditions for (\ref{eqthLIL.2})
are $\Sbep[X_1^2]<\infty$, $\Sbep[X_1]=\Sbep[-X_1]=0$ and  (\ref{eqthLIL.1}).
\end{remark}

The following corollary gives the compact law of the iterated logarithm.
\begin{corollary}\label{cor3.1}   Suppose that $X_1, X_2, \ldots$ are  independent  with
      $\lim\limits_{c\to \infty}  \Sbep\left[(X_1^2-c)^+\right]=0$, $\Sbep[X_1]=\Sbep[-X_1]=0$ and (\ref{eqthLIL.1}). Let $C \{x_n\} $ denote the cluster set of a sequence of $\{x_n\}$ in $\mathbb R$.  If $\Capc$ is continuous,
   then  we have
   \begin{description}
    \item[\rm (I)]
$$\Capc\left(C\Big\{\frac{S_n}{a_n}\Big\}=\Big[\liminf_{n\to \infty}\frac{S_n}{a_n},\;\;\limsup_{n\to \infty}\frac{S_n}{a_n}\Big]=
[-\overline{\sigma}, \;\;\overline{\sigma}]\right)=1; $$
    \item[\rm (II)]
    $$ \cCapc\left(\underline{\sigma}\le  \limsup_{n\to \infty}\frac{S_n}{a_n}\le \overline{\sigma}\right)=1 $$
    and
    $$ \cCapc\left(-\overline{\sigma}\le  \liminf_{n\to \infty}\frac{S_n}{a_n}\le -\underline{\sigma}\right)=1; $$
  \item[\rm (III)]
$$\cCapc\left([-\overline{\sigma}, \;\;\overline{\sigma}]\supset C\Big\{\frac{S_n}{a_n}\Big\}=\Big[\liminf_{n\to \infty}\frac{S_n}{a_n},\;\;\limsup_{n\to \infty}\frac{S_n}{a_n}\Big]\supset
[-\underline{\sigma}, \;\;\underline{\sigma}]\right)=1. $$
\end{description}
 If $\outCapc$ is continuous and $\Sbep$ is countably sub-additive,
   then we also have the conclusions (I), (II) and (III)   with $\cCapc$ being replaced $\outcCapc$.
\end{corollary}

\begin{remark} Note (I) and (III). It is interesting to know whether the cluster set $C\Big\{\frac{S_n}{a_n}\Big\}$ is the same  or not under $\Capc$ and $\cCapc$.
\end{remark}

\begin{remark}
 In the proof of Corollary~\ref{cor3.1}, the central limit theorems  under the Peng's framework  are used.
  Though a lot of results on the central limit theorems for related classical negatively dependent random variables can be found in literature
  (cf.  Newman \cite{Newman84}, Newman   and Wright \cite{NW81}, Su, Zhao and Wang \cite{SZW97}, Zhang \cite{Zh01b}  etc),
  it is difficult to establish a  central limit theorem  for non-independent random  variables under the sub-linear expectations.
  So, we have no version  of  Corollary~\ref{cor3.1} for negatively dependent random variables.
\end{remark}

\begin{remark} \label{remarkChenHu} Chen and Hu \cite{ChenHu14} tried to establish the law of the iterated logarithm for {\em bounded  }
 independent and identically distributed random variables. They have introduced a clever method to obtain the lower
bound $\underline{\sigma}$ of $\limsup\limits_{n\to \infty}\frac{S_n}{a_n}$.
 To establish the upper bound $\overline{\sigma}$,
 they have to assume the boundness of the random variables.
 Chen and Hu \cite{ChenHu14}  obtained   the exponential inequalities by taking $\varphi(x)=e^{\lambda x^2}$ in (\ref{eqCLT}). However,  $\varphi(x)=e^{\lambda x^2}$ is not a local Lipschitz function.  Similarly as in the study of  the moment convergence
  under the classical expectation/integral,
  when taking $\varphi(x)=e^{\lambda x^2}$ in (\ref{eqCLT}), the uniform integrability
  of $\exp\{\lambda(S_n/\sqrt{n})^2\}$ needs to be verified at first, which is not an easy work even in the case of the classical linear expectations. So there is a gap in Chen and Hu's proof.
\end{remark}

  \section{Proofs}\label{sectProof}
  \setcounter{equation}{0}

  We first show the exponential inequities, then the central limit theorem, and at last the law of the iterated logarithm.  The  H\"older's inequality under the sub-linear expectation will be used frequently in our proofs,
 which can be proved by the same may under the linear expectation due to the properties of
the monotonicity and sub-additivity, and the elementary inequality $|xy|\le \frac{1}{p}|x|^p+\frac{1}{q}|y|^q$, where $p,q>1$ are two real numbers satisfying $\frac{1}p+\frac{1}{q}=1$ (c.f, Proposition 16 of Denis,  Hu  and Peng \cite{DHP11}).
 \begin{lemma} ({\em H\"older's inequality}) Let $p,q>1$ be two real numbers satisfying $\frac{1}p+\frac{1}{q}=1$. Then for two random variables  $X,Y$  in $(\Omega, \mathscr{H}, \Sbep)$ we have
 $$ \Sbep[|XY|]\le \left(\Sbep[|X|^p]\right)^{\frac{1}{p}}  \left(\Sbep[|Y|^q]\right)^{\frac{1}{q}}$$
whenever $\Sbep[|X|^p]<\infty$, $\Sbep[|Y|^q]<\infty$.
 \end{lemma}

  We also   need
 the  Rosenthal type inequalities under $\Sbep$ which have been obtained  by Zhang \cite{Zh15b}.
 \begin{lemma}\label{thRIneq} ({\em Rosenthal's inequality})
(a) Let $\{X_1,\ldots, X_n\}$ be a sequence  of independent random variables
 in $(\Omega, \mathscr{H}, \Sbep)$ with $\Sbep[X_k]\le 0$, $k=1,\ldots, n$.
    Then
\begin{equation}\label{eqthRIneq.2}
\Sbep\left[\left|\max_{k\le n}(S_n- S_k)\right|^p\right]\le C_p\left\{ \sum_{k=1}^n \Sbep [|X_k|^p]+\left(\sum_{k=1}^n \Sbep [|X_k|^2]\right)^{p/2}\right\}, \;\; \text{ for }   p\ge 2.
\end{equation}

In particular,
\begin{equation}\label{eqthRIneq.3}
\Sbep\left[\left( S_n^+\right)^p\right]\le C_p\left\{ \sum_{k=1}^n \Sbep [|X_k|^p]+\left(\sum_{k=1}^n \Sbep [|X_k|^2]\right)^{p/2}\right\},
\;\; \text{ for }   p\ge 2.
\end{equation}
(b) Let $\{X_1,\ldots, X_n\}$ be a sequence  of negatively dependent random variables
 in $(\Omega, \mathscr{H}, \Sbep)$.
    Then
\begin{align}\label{eqthRIneq.4}
\Sbep\left[\max_{k\le n} \left|S_k\right|^p\right]\le  & C_p\left\{ \sum_{k=1}^n \Sbep [|X_k|^p]+\left(\sum_{k=1}^n \Sbep [|X_k|^2]\right)^{p/2} \right. \nonumber
 \\
& \qquad \left. +\left(\sum_{k=1}^n \big[\big(\cSbep [X_k]\big)^-+\big(\Sbep [X_k]\big)^+\big]\right)^{p}\right\}, \;\; \text{ for }   p\ge 2.
\end{align}
\end{lemma}

\bigskip

\subsection{Proofs of the exponential  inequalities}

{\bf Proof of     Theorem~\ref{th2}}.  Let $Y_k=X_k\wedge y$, $T_n=\sum_{k=1}^n Y_k$.  Then $X_k-Y_k=(X_k-y)^+\ge 0$ and
 $\Sbep[Y_k]\le \Sbep[X_k]\le 0$.
Note that $\varphi(x)=:e^{t(x\wedge y)}$ is a bounded non-decreasing function and belongs to $C_{l,Lip}(\mathbb R)$ since $0\le \varphi^{\prime}(x)\le t e^{ty}$ if $t>0$.
  It follows that for any $t>0$,
 $$
 \Capc\left(S_n\ge x\right)\le  \Capc\big(\max_{k\le n}X_k\ge y\big)+
 \Capc\left(T_n\ge x\right).
 $$
 and
\begin{align*}
\Capc\left(T_n\ge x\right)  \le
 e^{-tx}\Sbep[e^{t T_n}]
 \le
 e^{-tx}\prod_{k=1}^n \Sbep[e^{t Y_k}],
 \end{align*}
be the definition of the negative dependence. Note
$$e^{tY_k}=1+ tY_k+\frac{e^{tY_k}-1-t Y_k}{Y_k^2}Y_k^2\le 1 +tY_k+\frac{e^{ty}-1-t y}{y^2}Y_k^2.  $$
We have
$$ \Sbep[e^{t Y_k}]\le 1+\frac{e^{ty}-1-t y}{y^2} \Sbep[Y_k^2]\le\exp\left\{\frac{e^{ty}-1-t y}{y^2}\Sbep[X_k^2]\right\}. $$
Choosing $t=\frac{1}{y}\ln \big(1+\frac{xy}{B_n}\big)$ yields
\begin{align}\label{eqproofth2.1}
\Capc\left(T_n\ge x\right)  \le & e^{-tx}\exp\left\{\frac{e^{ty}-1-t y}{y^2}B_n\right\}\nonumber \\
=&\exp\left\{\frac{x}{y}-\frac{x}{y}\Big(\frac{B_n}{xy}+1\Big)\ln\Big(1+\frac{xy}{B_n}\Big)\right\}.
 \end{align}
 Applying the elementary inequality
 $$ \ln (1+t)\ge \frac{t}{1+t}+\frac{t^2}{2(1+t)^2}\big(1+\frac{2}{3} \ln (1+t)\big)$$
 yields
 $$ \Big(\frac{B_n}{xy}+1\Big)\ln\Big(1+\frac{xy}{B_n}\Big)
 \ge 1+\frac{xy}{2(xy+B_n)}\Big(1+\frac{2}{3}\ln\big(1+\frac{xy}{B_n}\big)\Big). $$
 (\ref{eqth2.1}) is proved.

 Next we show (b). If $xy\le \delta B_n$, then
 $$\frac{x^2}{2(xy+B_n)}\Big(1+\frac{2}{3}\ln\big(1+\frac{xy}{B_n}\big)\Big)\ge \frac{x^2}{2 B_n(1+\delta)}. $$
If  $xy\ge \delta B_n$,
then
 $$\frac{x^2}{2(xy+B_n)}\Big(1+\frac{2}{3}\ln\big(1+\frac{xy}{B_n}\big)\Big)\ge \frac{x}{2(1+1/\delta)y}. $$
 It follows that
\begin{equation}\label{eqproofth2.2} \Capc\left(T_n\ge x\right)\le \exp\left\{-\frac{x^2}{2 B_n (1+\delta)}\right\}
 +\exp\left\{-\frac{x}{2 (1+1/\delta)y}\right\}
 \end{equation}
 by (\ref{eqproofth2.1}). Let
 $$\beta(x)=\beta_p(x)=\frac{1}{x^p}\sum_{k=1}^n\Sbep[(X_k^+)^p], $$
 and choose
 $$ \rho=1\wedge \frac{1}{2(1+1/\delta)\delta \log (1/\beta(x))}, \;\; y=\rho\delta x. $$
 Then by (\ref{eqproofth2.2}),
 \begin{align*}
  &\Capc\big(S_n\ge (1+2\delta)x\big)\le  \Capc\big(T_n\ge  x\big)+\Capc\big(\sum_{i=1}^n(X_i-\rho\delta x)^+\ge 2\delta x\big) \\
  \le &\exp\left\{-\frac{x^2}{2 B_n (1+\delta)}\right\}+\beta(x)+\Capc\big(\max_{i\le n} X_i\ge \delta x\big)+\Capc\big(\sum_{i=1}^n(X_i-\rho\delta x)^+\wedge (\delta x) \ge 2\delta x\big).
  \end{align*}
  It is obvious that
  $$ \Capc\big(\max_{i\le n} X_i\ge \delta x\big)\le \delta^{-p} \beta(x). $$
  On the other hand, for $t>0$, 
  \begin{align*}
  &\Capc\big(\sum_{i=1}^n(X_i-\rho\delta x)^+\wedge (\delta x) \ge 2\delta x\big)=
 \Capc\left(\sum_{i=1}^n\left[\Big(\frac{X_i}{\delta x}-\rho\big)^+\wedge 1\right] \ge 2\right)\\
 \le & e^{-2t}\Sbep\exp\left\{t\sum_{i=1}^n\left[\Big(\frac{X_i}{\delta x}-\rho\big)^+\wedge 1\right]\right\}
 \le e^{-2t}\prod_{i=1}^n \Sbep\exp\left\{t\left[\Big(\frac{X_i}{\delta x}-\rho\big)^+\wedge 1\right]\right\}\\
 \le & e^{-2t}\prod_{i=1}^n\left[ 1+e^t \Capc(X_i\ge \rho\delta x)\right]
 \le \exp\left\{ -2t +e^t \sum_{i=1}^n \Capc(X_i\ge \rho\delta x)\right\}.
 \end{align*}
 Assume  $\beta(x)<1$. Suppose $\sum_{i=1}^n \Capc(X_i\ge \rho\delta x)<2$. Let $t=-\ln\frac{ \sum_{i=1}^n \Capc(X_i\ge \rho\delta x)}{2}$ (while, if $\sum_{i=1}^n \Capc(X_i\ge \rho\delta x)=0$, we let $t\to \infty$). We obtain
\begin{align}\label{eqproofth2.3}
&  \Capc\big(\sum_{i=1}^n(X_i-\rho\delta x)^+\wedge (\delta x) \ge 2\delta x\big)
\le   e^2 \left(\frac{1}{2}\sum_{i=1}^n \Capc(X_i\ge \rho\delta x)\right)^2\\
\le & e^2 \left(\frac{\beta(x)}{2(\delta \rho)^p}\right)^2
=e^22^{-2}\delta^{-2p} (2(\delta+1))^{2p}\beta^2(x)\left(\log\frac{1}{\beta(x)}\right)^{2p}
\le    C_p\delta^{-2p}\beta(x),\nonumber
\end{align}
where the last inequaltiy is due to the fact that $(\log 1/t)^{2p}\le C_p /t$ ($0<t<1$). When $\sum_{i=1}^n \Capc(X_i\ge \rho\delta x)\ge 2$, (\ref{eqproofth2.3}) is obvious.  So, we conclude  that
\begin{align*}
  \Capc\big(S_n\ge (1+2\delta)x\big)\le \exp\left\{-\frac{x^2}{2 B_n (1+\delta)}\right\}+C_p \delta^{-2p}\beta(x).
  \end{align*}
 If $\beta(x)\ge 1$, then the above inequality is obvious. Now letting $z= (1+2\delta)x $ and $\delta^{\prime}=(1+\delta)(1+2\delta)^2-1$ yields
 \begin{align*}
  \Capc\big(S_n\ge z\big)\le \exp\left\{-\frac{z^2}{2 B_n (1+\delta^{\prime})}\right\}+C_p (\delta^{\prime})^{-2p}\beta(z).
  \end{align*}
 (b) is proved.

Finally, we consider (c). Note that
$$ C_{\Capc}\big[(X^+)^p\big]= \int_0^{\infty}\Capc\big(X^p>x)dx= \int_0^{\infty}px^{p-1}\Capc\big(X>x)dx. $$
We put $y=x/r$, where $r=p>p/2$, in (\ref{eqproofth2.1}), then multiply both sides of this inequality by $px^{p-1}$. We find that
$$ px^{p-1}\Capc\left(S_n^+\ge x\right)\le  px^{p-1}\Capc\big(\max_{k\le n}X_k^+\ge \frac{x}{r}\big)
+pe^rx^{p-1}\left(1+\frac{x^2}{r B_n}\right)^{-r}. $$
By integrating on the positive half-line, we conclude (\ref{eqth2.3}).  $\Box$

\bigskip
{\bf Proof of    Corollary~\ref{cor2}}. Note, by (\ref{eq1.4}),
$$
 \cCapc\left(S_n\ge x\right)\le  \Capc\big(\max_{k\le n}X_k\ge x\big)+
 \cCapc\left(T_n\ge x\right), $$
  $$\cCapc\big(S_n\ge (1+2\delta)x\big)\le  \cCapc\big(T_n\ge  x\big)+\Capc\big(\sum_{i=1}^n(X_i-\rho\delta x)^+\ge 2\delta x\big),$$
$$\cSbep[e^{tY_k}]\le 1 +t\cSbep[Y_k]+\frac{e^{ty}-1-t y}{y^2}\Sbep[Y_k^2]\le   1 + \frac{e^{ty}-1-t y}{y^2}\Sbep[Y_k^2], $$
and $\cCapc(T_n\ge x)\le e^{-tx}\cSbep[e^{tT_n}]=  e^{-tx}\prod_{k=1}^n\cSbep[e^{tY_k}]$.
 The proof is similar to that of Theorem~\ref{th2}. $\Box$
\bigskip

For proving Theorem \ref{thLower}, we need the following  lemma on the $G$-normal distributed random variable, the proof of which can be found in Denis, Hu and Peng \cite{DHP11}.
\begin{lemma}\label{lemmaGnormal}
Suppose   $\xi\sim N(0, [\underline{\sigma}^2, \overline{\sigma}^2])$ under $\widetilde{\mathbb E}$. Let $P$ be a probability
measure and $\varphi$ be a bounded continuous function. If $\{B_t\}_{t\ge 0}$  is a $P$-Brownian motion,
then
$$ \widetilde{\mathbb E}\left[\varphi(\xi)\right]=\sup_{\theta_{\cdot}\in \Theta}\ep_P\left[\varphi\left(\int_0^1\theta_s dB_s\right)\right], $$
where
\begin{eqnarray*}
&\Theta=\left\{ \theta_{\cdot}:\theta_t \text{ is } \mathscr{F}_t\text{-adapted process such that }  \underline{\sigma}\le \theta_t\le \overline{\sigma}\right\},&\\
& \mathscr{F}_t=\sigma\{B_s:0\le s\le t\}\vee \mathscr{N}, \;\; \mathscr{N} \text{ is the collection of } P\text{-null subsets}. &
\end{eqnarray*}
\end{lemma}
\bigskip
\
{\bf Proof of      Theorem \ref{thLower}}.
It is important to note that the independence under $\Sbep$ is defined through continuous functions in $C_{l,Lip}$
 and the indicator function of an event is not continuous.  We need to modify the indicator function  by functions in $C_{l,Lip}$.
   For $t>0$, let
\begin{eqnarray*}
 N=:[n t^2/y_n^2],\;\; m=[y_n^2/t^2];\;\; r=\sqrt{n}y_n/(tm)
\end{eqnarray*}
and let $g_{\epsilon}$ be a function satisfying that its derivatives of each order are bounded, $g_{\epsilon}(x)=1$ if $x\ge 1$, $g_{\epsilon}(x)=0$ if $x\le 1-\epsilon$, and $0\le g_{\epsilon}(x) \le 1$ for all $x$,
where $0<\epsilon<1$.
Then
\begin{equation}\label{eqgfunction}
 g_{\epsilon}(\cdot) \in C_b^{\infty}(\mathbb R)\subset C_{l,Lip}(\mathbb R)\;  \text{ and }\; I\{x\ge 1\}\le g_{\epsilon}(x)\le I\{x>1-\epsilon\}.
 \end{equation}
Define a function
$ \phi(x)=  1-g_{1/2}\left(\frac{2 |x|}{\epsilon}\right). $
 Then
\begin{align*}
& \left\{ \left|  \frac{S_n}{y_n\sqrt{n} }-b\right|\le \epsilon\right\}
\supset   \left\{ b-\epsilon/2\le \frac{S_{Nm}}{y_n\sqrt{n} }\le b+ \epsilon/2\right\}
\bigcap \left\{ \left|\frac{S_n-S_{Nm}}{y_n\sqrt{n} }\right|\le \epsilon/2\right\}\\
& \;\;=    \left\{ tm(b-\epsilon/2)\le \frac{S_{Nm}}{r }\le tm(b+ \epsilon/2)\right\}
\bigcap \left\{ \left|\frac{S_n-S_{Nm}}{y_n\sqrt{n} }\right|\le \epsilon/2\right\}\\
&\;\; \supset   \bigcap_{i=1}^m  \left\{ bt-t\epsilon /2\le \frac{S_{Ni}-S_{N(i-1)}}{r }\le bt+ \epsilon t/2\right\}
\bigcap \left\{ \left|\frac{S_n-S_{Nm}}{y_n\sqrt{n} }\right|\le \epsilon/2\right\}\\
&\;\; =   \bigcap_{i=1}^m  \left\{  \Big| \frac{S_{Ni}-S_{N(i-1)}}{rt }-b\Big|\le   \epsilon  /2\right\}
\bigcap \left\{ \left|\frac{S_n-S_{Nm}}{y_n\sqrt{n} }\right|\le \epsilon/2\right\}
\end{align*}
It follows that
$$ I\left\{ \left|  \frac{S_n}{y_n\sqrt{n} }-b\right|\le \epsilon\right\}
\ge \prod_{i=1}^m \phi\left(\frac{S_{Ni}-S_{N(i-1)}}{rt }-b\right)
\left\{1-g_{1/2}\left(\frac{2}{\epsilon}\left|\frac{S_n-S_{Nm}}{y_n\sqrt{n} }\right|\right)\right\}. $$
Note that
 $\{S_{Ni}-S_{N(i-1)},i=1,\ldots,m, S_n-S_{Nm}\}$ are independent under $\Sbep$ (and $\cSbep$). By (\ref{eq1.3}),  we have
$$\cCapc\left( \left|  \frac{S_n}{y_n\sqrt{n} }-b\right|\le \epsilon\right)
\ge \left(\cSbep\left[ \phi\big(S_N/(rt)-b\big)\right]\right)^m
\left\{1-\Sbep\left[g_{1/2}\Big(\frac{2}{\epsilon} \frac{|S_{n-Nm}|}{y_n\sqrt{n} }\Big)\right]\right\}.
$$
Note
$$ 0\le \Sbep\left[g_{1/2}\Big(\frac{2}{\epsilon} \frac{|S_{n-Nm}|}{y_n\sqrt{n} }\Big)\right]\le
\Capc\left(\frac{|S_{n-Nm}|}{y_n\sqrt{n} }\ge \epsilon/4\right)\le \frac{16\overline{\sigma}^2}{\epsilon^2}\frac{n-Nm}{n y_n}\to 0.$$
By applying Theorem \ref{thCLT}, it follows that
\begin{align*}
& \liminf_{n\to \infty} y_n^{-2}\ln \cCapc\left( \left|  \frac{S_n}{y_n\sqrt{n} }-b\right|\le \epsilon\right) \\
\ge & \liminf_{t\to \infty}\liminf_{n\to \infty} t^{-2} m^{-1}\ln \cCapc\left( \left|  \frac{S_n}{y_n\sqrt{n} }-b\right|\le \epsilon\right) \\
\ge & \liminf_{t\to \infty}\liminf_{n\to \infty} t^{-2}  \ln \cSbep\left[ \phi\big(S_N/(rt)-b\big)\right] \\
=& \liminf_{t\to \infty} t^{-2}\ln \widetilde{\mathcal E}\left[ \phi\big(\xi/t-b\big)\right].
\end{align*}
Note $\xi/t\sim N(0,[\underline{\sigma}^2/t^2, \overline{\sigma}^2/t^2])$ under $\widetilde{\mathbb E}$. By Lemma 5   of Chen and Hu \cite{ChenHu14},
$$ \widetilde{\mathcal E}\left[ \phi\big(\xi/t-b\big)\right]
\ge \exp\left\{-\frac{1}{2}\Big(\frac{bt}{\underline{\sigma}}\Big)^2\right\}\widetilde{\mathcal E}\left[ \phi\big(\xi/t\big)\right]
\ge \exp\left\{-\frac{1}{2}\Big(\frac{bt}{\underline{\sigma}}\Big)^2\right\}\widetilde{\mathcal V}\Big(|\xi|\le \epsilon t/4\Big). $$
It follows that
$$\liminf_{t\to \infty} t^{-2}\ln \widetilde{\mathcal E}\left[ \phi\big(\xi/t-b\big)\right]\ge -\frac{1}{2}\Big(\frac{b}{\underline{\sigma}}\Big)^2.$$
 The proof of Theorem~\ref{thLower} (a)  is completed.

 For (b), with a similar argument we have
 $$\Capc\left( \left|  \frac{S_n}{y_n\sqrt{n} }-b\right|\le \epsilon\right)
\ge \left(\Sbep\left[ \phi\big(S_N/(rt)-b\big)\right]\right)^m
\left\{1-\cSbep\left[g_{1/2}\Big(\frac{2}{\epsilon} \frac{|S_{n-Nm}|}{y_n\sqrt{n} }\Big)\right]\right\},
$$
and then
$$
  \liminf_{n\to \infty} y_n^{-2}\ln \Capc\left( \left|  \frac{S_n}{y_n\sqrt{n} }-b\right|\le \epsilon\right)
\ge  \liminf_{t\to \infty} t^{-2}\ln \widetilde{\mathbb E}\left[ \phi\big(\xi/t-b\big)\right].
$$
From Lemma \ref{lemmaGnormal}, it follows that
\begin{align*}
\liminf_{t\to \infty} t^{-2}\ln \widetilde{\mathbb E}\left[ \phi\big(\xi/t-b\big)\right]
=&\liminf_{t\to \infty} t^{-2}\sup_{\theta_{\cdot}\in \Theta} \ln E_p\left[ \phi\big(\int_0^1\theta_s d B_s/t-b\big)\right] \\
\ge & \liminf_{t\to \infty} t^{-2}  \ln E_p\left[ \phi\big(  \overline{\sigma} B_1/t-b\big)\right]
\ge  -\frac{1}{2}\Big(\frac{b}{\overline{\sigma}}\Big)^2.
\end{align*}
The proof is completed. $\Box$.

\subsection{ Proofs of the central limit theorem}

For showing Theorem \ref{thCLT}, we let $Y_j=(-\sqrt{j})\vee \big(X_j\wedge \sqrt{j}\big)$, $T_n=\sum_{j=1}^n Y_j$.
 Suppose  that $\varphi$ is a bounded and (global)  Lipschitz continuous function.  We first show that
\begin{equation}\label{eqproofthCLT1.1}
 \lim_{n\to \infty}\Sbep\left[\varphi\left(\frac{T_n}{\sqrt{n}}\right)\right]=\widetilde{\mathbb E}[\varphi(\xi)].
 \end{equation}
We use the  argument of Peng \cite{Peng08b} for a bounded and  Lipschitz continuous function, which is a version of the Stein method under the sub-linear expectations.  Here we only give the difference.
The main difference is that   $Y_j$s are not identically distributed and
$\Sbep[Y_j]$, $\Sbep[-Y_j]$   are no zeros.

First, we have the following facts.
\begin{description}
\item[\rm (F1)] $\Sbep[(X_1^2-j)^+]\to 0$ as $j\to \infty$, and
\begin{equation}\label{eqCLTfact1}
 \frac{\sum_{j=1}^n \Sbep[|X_j-Y_j|]}{\sqrt{n}}\to 0 \;\text{ as } n\to \infty.
\end{equation}
\item[\rm (F2)]
\begin{equation}\label{eqCLTfact2} \frac{\sum_{j=1}^n \Sbep\big[|Y_j|^{2+\alpha}\big]}{n^{1+\alpha/2}}\to 0
 \;\text{ as } n\to \infty, \;\; \forall \alpha>0.
 \end{equation}
\item[\rm (F3)] For any $p\ge 2$,
\begin{equation}\label{eqCLTfact3}\Sbep[|T_n|^p]\le C_p n^{p/2}.
\end{equation}
\end{description}
In fact, for (F1), note
$$ \Sbep|X_j-Y_j|\le \Sbep\left[(|X_1|-\sqrt{j})^+\right]
\le j^{-1/2}\Sbep\left[(X_1^2-j)^+\right].$$
(F1) is obvious.

For (F2), note that $\Sbep\big[|Y_j|^{2+\alpha}\big]\le c^{2+\alpha} +j^{\alpha/2}\Sbep[(X_1^2-c)^+]$ for any $c>1$. So, (\ref{eqCLTfact2}) is true.

For (F3), by the Rosenthal inequality (\ref{eqthRIneq.4}) and the fact (\ref{eqCLTfact1}) we have
\begin{align*}
 \Sbep\left[|T_n|^p\right]\le &
  C_p\sum_{j=1}^n \Sbep[|Y_j|^p]+C_p\Big(\sum_{j=1}^n \Sbep[|Y_j|^2]\Big)^{p/2}
+ C_p \Big( \sum_{j=1}^n\big[(\Sbep Y_j)^++(\Sbep[- Y_j])^+ \Big)^p \\
\le & C_p n^{p/2-1} \sum_{j=1}^n   \Sbep[ X_j^2]+C_p\Big( \sum_{j=1}^n \Sbep[|Y_j|^2]\Big)^{p/2}
+C_p\Big( \sum_{j=1}^n\big[ \Sbep [|X_j-Y_j|]  \Big)^p
\le   C_p n^{p/2}.
\end{align*}
  And so, (\ref{eqCLTfact3}) is true.

Now,  for a small but fixed $h > 0$, let $V (t, x)$ be the unique viscosity solution of the following equation,
$$ \partial_t V + G( \partial_{xx}^2 V)=0,\;\;  (t, x) \in [0,1+ h] \times \mathbb R, \; V|_{t=1+h} = \varphi(x), $$
where $G(\alpha)=\frac{1}{2}\big(\overline{\sigma}^2 \alpha^+-\underline{\sigma}^2\alpha^-\big)$. Then by
the interior regularity of $V$,
$$ \|V\|_{C^{1+\alpha/2,2+\alpha}([0,1]\times R)} < \infty, \text{ for some } \alpha\in (0, 1). $$
According to the definition of $G$-normal distribution, we have $V(t,x)=\widetilde{\mathbb E}\big[\varphi(x+\sqrt{1+h-t}\xi)\big]$. In particular,
$$ V(h,0)=\widetilde{\mathbb E}\big[\varphi( \xi)\big], \;\; V(1+h, x)=\varphi(x). $$
It is obvious that, if  $\varphi(\cdot)$ is a global Lipschitz function, i.e.,
$|\varphi(x)-\varphi(y)|\le C|x-y|$, then $|V(t,x)-V(t,y)|<C|x-y|$ and $|V(t,x)-V(s,x)|\le C\widetilde{\mathbb E}[|\xi|] |t-s|^{1/2} $.
So, $|V(1+h,x)-V(1,x)|\le C\widetilde{\mathbb E}[|\xi|] \sqrt{h}$ and $|V(h,0)-V(0,0)|\le C\widetilde{\mathbb E}[|\xi|] \sqrt{h}$. Let $\delta=\frac{1}{n}$, $T_0=0$.  Following the proof of Lemma 5.4 of Peng \cite{Peng08b}, it is sufficient to show that
\begin{equation}\label{eqconV}
 \lim_{n\to \infty} \Sbep[ V(1,\sqrt{\delta}T_n)]=V(0,0).
 \end{equation}
Applying the Taylor's expansion yields
\begin{align*}
 & V(1,\sqrt{\delta}T_n)-V(0,0)\\
 =&\sum_{i=0}^{n-1} \left\{[ V((i+1)\delta, \sqrt{\delta}T_{i+1})-V(i\delta, \sqrt{\delta}T_{i+1})]
  +[ V(i\delta, \sqrt{\delta}T_{i+1})-V(i\delta, \sqrt{\delta}T_i)]\right\}\\
 =& \sum_{i=0}^{n-1} \left\{I_{\delta}^i+J_{\delta}^i\right\},
\end{align*}
with
$  |I_{\delta}^i|\le C\delta^{1+\alpha/2}\big(1+|Y_{i+1}|^{\alpha}+|Y_{i+1}|^{2+\alpha}\big), $
\begin{align*}
 J_{\delta}^i =& \partial_t V(i\delta,\sqrt{\delta} T_i)\delta +\frac{1}{2} \partial_{xx}^2 V(i\delta,\sqrt{\delta} T_i)Y_{i+1}^2 \delta+
\partial_x V(i\delta,\sqrt{\delta} T_i) Y_{i+1}\sqrt{\delta} \\
=& \left(\partial_t V(i\delta,\sqrt{\delta} T_i)\delta +\frac{1}{2} \partial_{xx}^2 V(i\delta,\sqrt{\delta} T_i)X_{i+1}^2 \delta+
\partial_x V(i\delta,\sqrt{\delta} T_i) X_{i+1}\sqrt{\delta}\right) \\
&+\left( \frac{1}{2} \partial_{xx}^2 V(i\delta,\sqrt{\delta} T_i)\big(Y_{i+1}^2-X_{i+1}^2\big) \delta+
\partial_x V(i\delta,\sqrt{\delta} T_i)\big( Y_{i+1}-X_{i+1}\big)\sqrt{\delta} \right)\\
=:&J_{\delta,1}^i+J_{\delta,2}^i,
\end{align*}
where $C$ is a constant.  And so
\begin{align*}
&\Sbep[\sum_{i=0}^{n-1}J_{\delta,1}^i]-\sum_{i=0}^{n-1}\Big\{\Sbep[|J_{\delta,2}^i|]+\Sbep[|I_{\delta}^i|]\Big\}\\
\le &\Sbep[V(1,\sqrt{\delta}T_n)]-V(0,0)\le \Sbep[\sum_{i=0}^{n-1}J_{\delta,1}^i]+\sum_{i=0}^{n-1}\Big\{\Sbep[|J_{\delta,2}^i|]+\Sbep[|I_{\delta}^i|]\Big\}.
\end{align*}
By noting the fact  (F2), we have
$$\sum_{i=0}^{n-1}\Sbep[|I_{\delta}^i|]\le C \Big(\frac{1}{n}\Big)^{1+\alpha/2}\sum_{i=1}^{n-1}\Big(1+\Sbep[|Y_i|^{\alpha}]+\Sbep[|Y_i|^{2+\alpha}]\Big)\to 0. $$
For $J_{\delta,1}^i$, note $\Sbep[X_{i+1}^2]=\overline{\sigma}^2$, $\cSbep[X_{i+1}^2]=\underline{\sigma}^2$,
$\Sbep[X_{i+1} ]=\cSbep[X_{i+1} ]=0$. It follows that
$$ \Sbep\left[J_{\delta,1}^i\big|X_1,\ldots, X_i\right]=
\big[\partial_t V(i\delta,\sqrt{\delta} T_i)+G\big(\partial^2_{xx}   V(i\delta,\sqrt{\delta} T_i)\big)\big]\delta=0. $$
It follows that $\Sbep\left[\sum_{i=0}^{n-1}J_{\delta,1}^i\right]=\Sbep\left[\sum_{i=0}^{n-2}J_{\delta,1}^i\right]=\ldots=0$.
For $J_{\delta,2}^i$,  since $\partial_{xx}V$ is uniformly $\alpha$-H\"older continuous in $x$ and
$\alpha/2$-H\"older continuous in $t$ on $[0,1]\times R$, we have
\begin{align*}
 &\Sbep\left[\big|\partial_{xx}^2 V(i\delta,\sqrt{\delta} T_i)\big|\right]\le \big|\partial_{xx}^2 V(0,0)\big|+ \Sbep\left[\big|\partial_{xx}^2 V(i\delta,\sqrt{\delta} T_i)-\partial_{xx}^2 V(0,0)\big|\right]\\
 \le & C\big[1+(i\delta)^{\alpha/2}+\Sbep[|\sqrt{\delta} T_i|^{\alpha}\big]
 \le   C\big[1+(i\delta)^{\alpha/2}+\big(\Sbep[|\sqrt{\delta} T_i|^2\big)^{\alpha/2}\big]\le C
 \end{align*}
by the H\"older inequality and Fact 3.  Similarly, $\Sbep\left[\big|\partial_{x} V(i\delta,\sqrt{\delta} T_i)\big|\right]\le C$.
It follows that
\begin{align*}
 \sum_{i=0}^{n-1}\Sbep[|J_{\delta,2}^i|]
\le &  \sum_{i=0}^{n-1}\left\{\frac{1}{2}\Sbep\left[\big|\partial_{xx}^2 V(i\delta,\sqrt{\delta} T_i)\big|\right]\Sbep[|X_{i+1}^2-Y_{i+1}^2|]\delta
\right.\\
& \left.+\Sbep\left[\big|\partial_{x} V(i\delta,\sqrt{\delta} T_i)\big|\right]\Sbep[|X_{i+1}-Y_{i+1}|]\sqrt{\delta}\right\}\\
\le & C\frac{1}{n}\sum_{i=1}^n \Sbep[(X_1^2-j)^+]+\frac{1}{\sqrt{n}}\sum_{i=1}^n\Sbep[|X_i-Y_i|]\to 0,
\end{align*}
by the independence and  Fact 1. (\ref{eqconV}) is now proved and hence (\ref{eqproofthCLT1.1}) follows.
 Finally, by the Lipschitz continuity of  $\varphi$, we have
\begin{align*}
  \left|\Sbep\left[\varphi\left(\frac{S_n}{\sqrt{n}}\right)\right]-\Sbep\left[\varphi\left(\frac{T_n}{\sqrt{n}}\right)\right]\right|
\le C\frac{\sum_{j=1}^n\Sbep[|X_j-Y_j|]}{\sqrt{n}} \to 0,
\end{align*}
by Fact (F1). So, for a  bounded and  Lipschitz continuous function $\varphi$, (\ref{eqCLT}) is verified.

If $\varphi$ is a bounded and  uniformly continuous function,
 we   define a function $\varphi_{\delta}$ as a convolution of $\varphi$ and the density of a normal distribution
$N(0,\delta)$, i.e.,
$$ \varphi_{\delta}= \varphi\ast \psi_{\delta}, \;\; \text{with} \;
\psi_{\delta}(x)=\frac{1}{\sqrt{2\pi\delta}}\exp\left\{-\frac{x^2}{2\delta}\right\}. $$
 Then $|\varphi_{\delta}^{\prime}(x)|\le \sup_x|\varphi(x)|\delta^{-1/2}$ and $\sup_x|\varphi_{\delta}(x)-\varphi(x)|\to 0$ as $\delta\to 0$.
 As proved, (\ref{eqCLT}) holds for each $\varphi_{\delta}$. So, it holds for $\varphi$.

Finally, suppose that $p\ge 2$, $\Sbep[|X_1|^p]<\infty$,  $\lim\limits_{c\to \infty} \Sbep[(|X_1|^2-c)^+]=0$,  and $\varphi$ is a continuous function satisfying $\varphi(x)\le C( 1+|x|^p)$. Give a number $N>1$.
Define $\varphi_1(x)=\varphi\big((-N)\vee (x\wedge N)\big)$ and $ \varphi_2(x)=\varphi(x) -\varphi_1(x)$.
 Then $\varphi_1$ is a bounded and uniformly continuous function and
$$ |\varphi_2(x)| \le 4C|x|^pI\{|x|>N\}\le 4C(2|x|^p-N)^+= 8C(|x|^p-N/2)^+. $$
So
\begin{align*}
 &\left|\Sbep\left[\varphi\left(\frac{S_n}{\sqrt{n}}\right)\right]-\widetilde{\mathbb E}[\varphi(\xi)]\right|\\
\le & \left|\Sbep\left[\varphi_1\left(\frac{S_n}{\sqrt{n}}\right)\right]-\widetilde{\mathbb E}[\varphi_1(\xi)]\right|
  + 8C\Sbep\left[ \left(  |S_n/\sqrt{n}|^p -N/2\right)^+\right]
  + 8C\Sbep\left[ \left(  |\xi|^p -N/2\right)^+\right].
\end{align*}
Hence, it is sufficient to show that
\begin{equation}\label{eqproofthCLT1.2}\lim_{N\to \infty}\limsup_{n\to \infty}\Sbep\left[ \left(  |S_n/\sqrt{n}|^p -N\right)^+\right]=0.
\end{equation}
Let $\widehat{Y}_j=X_j-Y_j$, , $\widehat{S}_n=\sum_{j=1}^n (\widehat{Y}_j-\Sbep[\widehat{Y}_j])$. Then
$$ S_n^+\le  T_n^+ +\widehat{S}_n^+ +\sum_{j=1}^n \Sbep[|\widehat{Y}_j|],$$
$$ \left( \left|\frac{S_n^+}{\sqrt{n}}\right|^p -N\right)^+\le
\left(3^{p-1}\left|\frac{T_n^+}{\sqrt{n}}\right|^p -N\right)^+ +3^{p-1}\left|\frac{\widehat{S}_n^+}{\sqrt{n}}\right|^p +3^{p-1} \left(\sum_{j=1}^n \frac{\Sbep[|\widehat{Y}_j|]}{\sqrt{n}}\right)^p. $$
We have shown that
$$\sum_{j=1}^n \frac{\Sbep[|\widehat{Y}_j|]}{\sqrt{n}}= \sum_{j=1}^n  \frac{\Sbep[|X_j-Y_j|]}{\sqrt{n}}\to 0 $$
by Fact (F1), and
$$\Sbep\left[\left(3^{p-1}\left|\frac{T_n^+}{\sqrt{n}}\right|^p -N\right)^+ \right]\le N^{-1} 3^{2p-2}\Sbep\left[\left|\frac{T_n}{\sqrt{n}}\right|^{2p}\right]\le N^{-1}C_p $$
by Fact (F3). Applying (\ref{eqthRIneq.3})   yields
\begin{align*}
  &\Sbep\left[|\widehat{S}_n^+/\sqrt{n}|^p\right]
\le
C_p  n^{ -p/2} \sum_{j=1}^n \Sbep[|\widehat{Y}_j|^p] +C_p  \Big(n^{ -1} \sum_{j=1}^n \Sbep[|\widehat{Y}_j|^2]\Big)^{p/2}\\
\le & C_p n^{-p/2} \sum_{j=1}^n\Sbep\left[\big((|X_1|-j^{1/2})^+\big)^p\right]+C_p  \Big(n^{ -1} \sum_{j=1}^n \Sbep\left[\big((|X_1|-j^{1/2})^+\big)^2\right]\Big)^{p/2}.
\end{align*}
When $p=2$,
$$\Sbep\left[|\widehat{S}_n^+/\sqrt{n}|^p\right]\le Cn^{ -1} \sum_{j=1}^n \Sbep\left[\big((|X_1|-j^{1/2})^+\big)^2\right]\to 0. $$
When $p>2$,
$$\Sbep\left[|\widehat{S}_n^+/\sqrt{n}|^p\right]
\le C_p n^{-p/2+1}  \Sbep\left[ |X_1|^p\right]+C_p  \Big(n^{ -1} \sum_{j=1}^n \Sbep\left[\big((|X_1|-j^{1/2})^+\big)^2\right]\Big)^{p/2}\to 0. $$
By combining the above arguments, it follows that
$$\lim_{N\to \infty}\limsup_{n\to \infty}\Sbep\left[ \left(  |S_n^+/\sqrt{n}|^p -N\right)^+\right]=0. $$
Similarly,
$$\lim_{N\to \infty}\limsup_{n\to \infty}\Sbep\left[ \left(  |S_n^-/\sqrt{n}|^p -N\right)^+\right]=0. $$
Hence (\ref{eqproofthCLT1.2}) is proved and the proof is now completed. $\Box$.

\bigskip

If consider $Y_j=(-j)\vee (X_j\wedge j)$ and  the following equation, instead,
$$ \partial_t V + \overline{G}( \partial_{x} V)=0, $$
where $\overline{G}(\alpha)= \overline{\mu} \alpha^+-\underline{\mu}\alpha^- $,
 one can prove the following weak law of large numbers. The  proof is similar to that of Theorem~\ref{thCLT} and so omitted.
 \begin{corollary}\label{thWLLN} (WLLN) Suppose that $\{X_n; n\ge 1\}$ is a sequence of   independent  and identically distributed random variables with
   $\Sbep[X_1]=\overline{\mu}$, $\cSbep[X_1]=\underline{\mu}$ and
   $\lim\limits_{c\to \infty} \Sbep\left[(|X_1|-c)^+\right]=0$.  Then for any   continuous function $\varphi$ satisfying  $|\varphi(x)|\le C (1+|x|)$,
  \begin{equation} \label{eqWLLN} \lim_{n\to \infty}\Sbep\left[\varphi\left(\frac{S_n}{n}\right)\right]
  =\sup_{\underline{\mu}\le x\le \overline{u}} \varphi(x).
  \end{equation}
    Further, if $p\ge 1$ and $ \Sbep[|X_1|^p]<\infty$,
     then (\ref{eqWLLN}) holds for any   continuous function $\varphi$ satisfying  $|\varphi(x)|\le C (1+|x|^p)$.
 \end{corollary}

  \subsection{Proofs of the laws of the iterated logarithm}

For proving the law of the iterated  logarithm, we need more properties of the sub-linear expectations and capacities.  We define an extension of $\Sbep$ on the space
 of all random variables by
$$ \extSbep[X]=\inf\{ \Sbep[Y]: X\le Y, Y\in \mathscr{H}\}. $$
Then $\extSbep$ is a sub-linear expectation on the space of all random variables, and
$$ \extSbep[X]=\Sbep [X] \;\; \forall X\in \mathscr{H}, \;\; \Capc(A)=\extSbep[I_A]\;\; \forall A\in \mathcal F. $$

\begin{lemma} \label{lem2} Suppose $X\in \mathscr{H}$.
\begin{description}
  \item[\rm (i)]
  Then for any $\delta>0$,
$$ \sum_{n=1}^{\infty} \Capc\big(|X|\ge \delta a_n\big)<\infty \;\; \Longleftrightarrow C_{\Capc}\left[\frac{X^2}{\log\log|X|}\right]<\infty.
$$
 \item[\rm (ii)]
   If $C_{\Capc}\left[\frac{X^2}{\log\log|X|}\right]<\infty$, then for any $\delta>0$ and $p>2$,
$$ \sum_{n=1}^{\infty} \frac{\Sbep\big[\big(|X|\wedge (\delta a_n)\big)^p\big]}{a_n^p}<\infty. $$
\item[\rm (iii)] For any $0\le b<c<\infty$,
\begin{equation}\label{eqlem2.3} \Sbep\big[|X|\wedge c\big]\le \int_0^c\Capc(|X|>x) dx), \;\;
\extSbep\big[|X|I\{b\le |X|\le  c\}\big]\le \int_b^c\Capc(|X|>x) dx.
\end{equation}
If $\lim\limits_{c\to +\infty}\Sbep[(|X|-c)^+]=0$ or  $\Sbep$ is countably sub-additive,
then $\Sbep[|X|]\le C_{\Capc}(|X|)$.
\end{description}
\end{lemma}
{\bf Proof}. (i) It is sufficient to note that
$$\left\{c_1\frac{|X|^2}{\log\log |X|}>n\right\}\subset \Big\{|X|\ge \delta a_n\Big\} \subset\left\{c_2\frac{|X|^2}{\log\log |X|}>n\right\}, \; n\ge n_0  $$
for some $0<c_1<c_2$ and $n_0$, and
$$ \sum_{n=1}^{\infty} \Capc\big(|Y|>n\big)<\infty \Longleftrightarrow \int_0^{\infty} \Capc\big(|Y|>x\big)dx. $$

(iii) Suppose that $Y$ is a bounded random variable. Choose  $n>2$ such that $|Y|<n$. Then
\begin{align*}
|Y|  = \sum_{i=1}^n|Y|I\{i-1<|Y|\le i\}
\le & \sum_{i=1}^n i\big(I\{|Y|> i-1\}-I\{|Y|> i\}\big) \\
\le & 1+\sum_{i=1}^nI\{|Y|> i\}.
\end{align*}
It follows that
\begin{equation}\label{eqprooflem2.1} \extSbep\big[ |Y| \big]\le 1+\sum_{i=1}^n\Capc(|Y|> i)\le 1+\int_0^{\infty} \Capc(|Y|> x)dx\le 1+  C_{\Capc}[|Y|]
\end{equation}
by the (finite) sub-additivity of $\extSbep$.
By considering $|Y|/\epsilon$ instead of $|Y|$, we have
$$ \extSbep\left[\frac{|Y|}{\epsilon}\right] \le 1+  C_{\Capc}\left[\frac{|Y|}{\epsilon}\right]=1+\frac{1}{\epsilon} C_{\Capc}[|Y|]. $$
That is
$ \extSbep\big[|Y|\big] \le \epsilon +  C_{\Capc}[|Y|].$
Taking $\epsilon\to 0$ yields
$ \extSbep\big[|Y|\big] \le   C_{\Capc}[|Y|]. $
Now, taking $Y=|X|\wedge c$ and $Y=|X|I\{b\le |X|\le c\}$ completes the proof of (\ref{eqlem2.3}).

If $\lim\limits_{c\to +\infty}\Sbep[(|X|-c)^+]=0$, taking $c\to \infty$ yields
$\Sbep\big[|X|\big]=\lim\limits_{c\to \infty}\Sbep\big[|X|\wedge c\big]\le C_{\Capc}(|X|)$.
 If $\Sbep$ is countably sub-additive, then so is $\extSbep$. Taking $n=\infty$ in (\ref{eqprooflem2.1})  and using the countable sub-additivity yields the same
conclusion for any $Y$.

(ii) Let $ f(x)$ be the inverse function of $\sqrt{2x\log\log x}$. By (c),
\begin{align*}
  \Sbep\big[\big(|X|\wedge (\delta a_n)\big)^p\big]\le  &\int_0^{(\delta a_n)^p}\Capc\left(|X|^p>x\right)dx
=p\int_0^{\delta a_n}x^{p-1}\Capc\left(|X|>x\right)dx \\
\le & 2p\int_0^{\delta n} (2\log\log y)^{\frac{p}{2}} y^{\frac{p}{2}-1} \Capc\left(f(|X|)>y\right)dy.
\end{align*}
It follows that
\begin{align*}
&\sum_{n=16}^{\infty} \frac{\Sbep\big[\big(|X|\wedge (\delta a_n)\big)^p\big]}{a_n^p}
\le  2p\sum_{n=16}^{\infty}  a_n^{-p}\int_0^{\delta n} (2\log\log y)^{\frac{p}{2}} y^{\frac{p}{2}-1} \Capc\left(f(|X|)>y\right)dy\\
\le &  4p\int_{16}^{\infty}(2x\log\log x)^{-p/2} \int_0^{\delta x} (2\log\log y)^{\frac{p}{2}} y^{\frac{p}{2}-1} \Capc\left(f(|X|)>y\right)dydx \\
\le &  4p\int_{0}^{\infty} (2\log\log y)^{\frac{p}{2}} y^{\frac{p}{2}-1} \Capc\left(f(|X|)>y\right)dy\int_{y/\delta}  (2x\log\log x)^{-p/2}dx \\
\le &c_{\delta}\int_{0}^{\infty}  \Capc\left(f(|X|)>y\right)dy=c_{\delta}C_{\Capc}\big[f(|X|)\big]
\le c_{\delta}C_{\Capc}\left[\frac{|X|^2}{\log\log |X|}\right]<\infty. \;\; \Box
\end{align*}

\bigskip
{\bf Proof of Theorem \ref{thLIL}.} The proof will be completed via three steps. The first step is to show (a).
 The second step is show that (\ref{eqthLIL.3}) implies (\ref{eqthLIL.1}). In the last step, we show that the means are zeros.

 \noindent{\bf Step 1}. We show (a).  Without loss of generality, we assume $\Sbep[X_1^2]=1$.  We will show that
\begin{equation}\label{eqLIL1.10}
\limsup_{n\to \infty} \frac{S_n}{a_n}\le 1 \;\; a.s. \; \outCapc, \;\; i.e.\;\; \outCapc\left(\Big\{\frac{S_n}{a_n}>1\Big\}\; i.o.\right)=0.
\end{equation}
 For given  $0<\epsilon<1$, let $b_n=\frac{\epsilon}{40} \sqrt{n/\log\log n}$, $a_n=\sqrt{2n\log\log n}$.
Let $Y_k=(-b_k)\vee (X_k\wedge b_k)$, $\widehat{Y}_k=X_k-Y_k$.   By the countable  sub-additivity of $\outCapc$,
(\ref{eqLIL1.10}) will follow if we have shown
that
\begin{equation}\label{eqLIL1.11}
\frac{\sum_{k=1}^n \widehat{Y}_k^+}{a_n}\to 0 \;\; a.s. \outCapc\;  \text{ and } \; \frac{\sum_{k=1}^n \widehat{Y}_k^-}{a_n}\to 0
\;\; a.s. \outCapc
\end{equation}
and
\begin{equation}\label{eqLIL1.12}
\limsup_{n\to \infty} \frac{\sum_{k=1}^n Y_k}{a_n}\le (1+\epsilon)^2\;\; a.s. \outCapc.
\end{equation}

We first show (\ref{eqLIL1.11}). We only consider $\widehat{Y}_k^+$ since the proof for $\widehat{Y}_k^-$ is similar.
Let $g_{\epsilon}(\cdot)$ be a smooth  function satisfying  (\ref{eqgfunction}). By (\ref{eq1.3}),
\begin{align*}
 \sum_{k=1}^{\infty}\Capc\big( \widehat{Y}_k^+> a_k\big)\le & \sum_{k=1}^{\infty}\Capc\big( X_k> a_k+b_k\big)
\le   \sum_{k=1}^{\infty}\Capc\big( X_k>  a_k\big)\\
\le & \sum_{k=1}^{\infty}\Sbep\left[g_{1/2}\big(X_k/a_k\big)\right]
= \sum_{k=1}^{\infty}\Sbep\left[g_{1/2}\big(X_1/a_k\big)\right]\\
\le &\sum_{k=1}^{\infty}\Capc\big( X_1> a_k/2\big)\le c C_{\Capc}\Big[\frac{|X_1|^2}{\log\log |X_1|}\Big]<\infty.
\end{align*}
Hence $\outCapc\big(\{\widehat{Y}_k^+>a_k\}\;i.o.\big)=0$ by the Borel-Cantelli lemma.
So, it is sufficient to show that
\begin{equation}\label{eqLIL1.13} \frac{\sum_{k=1}^n (\widehat{Y}_k^+)\wedge a_k}{a_n}\to 0 \;\; a.s. \outCapc.
\end{equation}
Note that the random variables $(\widehat{Y}_k^+)\wedge a_k$'s are non-negative. It is sufficient to show that
$$ \frac{\sum_{k=2^n+1}^{2^{n+1}} (\widehat{Y}_k^+)\wedge a_k}{a_{2^{n+1}}}\to 0 \;\; a.s.\; \outCapc. $$
Let $Z_k=(\widehat{Y}_k^+)\wedge a_k-\Sbep[(\widehat{Y}_k^+)\wedge a_k]$. Then for any $\delta>0$,
\begin{align*}
 \Capc\Big(\frac{\sum_{k=2^n+1}^{2^{n+1}} Z_k}{a_{2^{n+1}}}\ge \delta\Big)\le
 C\sum_{k=2^n+1}^{2^{n+1}}\frac{ \Sbep[|Z_k|^p]}{a_k^{p}}+\exp\left\{-\frac{\delta^2(a_{2^{n+1}})^2} {4\sum_{k=2^n+1}^{2^{n+1}}\Sbep Z_k^2}\right\}
\end{align*}
by (\ref{eqth2.2}). Note $\Sbep Z_k^2\le 4\Sbep\big[\big((X_1-b_k)^+\big)^2\big]\to 0$.  By Lemma \ref{lem2} (ii),  we have
\begin{align*}
\sum_n \Capc\Big(\frac{\sum_{k=2^n+1}^{2^{n+1}} Z_k}{a_{2^{n+1}}}\ge \delta\Big)\le
 C\sum_{k=1}^{\infty} \frac{\Sbep[(|X_k|\wedge a_k)^p]}{a_k^{p}}+C\sum_{n=1}^{\infty} \exp\big\{-2\log\log 2^{n+1}\big\}<\infty.
\end{align*}
Hence by the countable sub-additivity of $\outCapc$ and the Borel-Cantelli Lemma, we have
$$\outCapc\left(\Big\{\frac{\sum_{k=2^n+1}^{2^{n+1}} Z_k}{a_{2^{n+1}}}>\delta\Big\}\; i.o. \right)=0,\;\; \forall \delta>0. $$
On the other hand,
$$ \frac{\sum_{k=2^n+1}^{2^{n+1}} \Sbep\big[(\widehat{Y}_k^+)\wedge a_k\big]}{a_{2^{n+1}}}
\le C\frac{\sum_{k=2^n+1}^{2^{n+1}}\frac{(\log\log k)^{1/2}}{k^{1/2}} \extSbep[X_1^2I\{|X_1|\ge b_k\}]}{ a_{2^{n+1}}}\to 0.
$$
Hence, (\ref{eqLIL1.13}) is proved. So (\ref{eqLIL1.11}) holds.

Next, we show (\ref{eqLIL1.12}).
Let $n_k=[e^{k^{1-\alpha}}]$, where $0<\alpha<\frac{\epsilon}{1+\epsilon}$. Then $n_{k+1}/n_k\to 1$ and $\frac{n_{k+1}-n_k}{n_k}\approx \frac{C}{k^{\alpha}}$. For $n_k<n\le n_{k+1}$,
we have
\begin{align*}
\sum_{i=1}^n Y_k\le & \sum_{i=1}^{n_k} Y_i+\max_{n_k<n\le n_{k+1}}\Big(\sum_{i=n_k+1}^n Y_i\Big)\\
\le & \sum_{i=1}^{n_k} (Y_i-\Sbep[Y_i])+\max_{n_k<n\le n_{k+1}}\Big(\sum_{i=n_k+1}^{n} Y_i \Big)
  +\sum_{i=1}^{n_k} \Big| \Sbep[Y_i] \Big| \\
=:& I_k+II_k+III_k.
\end{align*}
For the third term, from the fact that $\Sbep[X_i]=\Sbep[-X_i]=0$ it follows that
$\big| \Sbep[Y_i] \big|=\big| \Sbep[Y_i] -\Sbep[X_i]\big|\le  \Sbep[(|X_1|-b_i)^+]\le b_i^{-1}\Sbep\big[(X_1^2- b_i^2)^+\big]$ and
$\big| \Sbep[-Y_i] \big|=\big| \Sbep[-Y_i] -\Sbep[-X_i]\big|\le \Sbep[(|X_1|-b_i)^+]\le b_i^{-1}\Sbep\big[(X_1^2-b_i^2)^+\big]$.
Hence
\begin{equation}\label{eqLIL1.15}
\frac{III_k}{a_{n_k}}=o(1)\frac{\sum_{i=1}^{n_k} \frac{(\log\log i)^{1/2}}{i^{1/2}}}{(n_k\log\log n_k)^{1/2}}\to 0.
\end{equation}

For the second term, by applying the Rosenthal inequality (\ref{eqthRIneq.4})  we have
\begin{align*}
\Capc\left(II_k\ge  \delta a_{n_k}\right)\le & c \frac{\sum_{i=n_k+1}^{n_{k+1}} \Sbep[|Y_i|^p]}{a_{n_k}^p}+
 c\left( \frac{\sum_{i=n_k+1}^{n_{k+1}} \Sbep[|Y_i|^2]}{a_{n_k}^2}\right)^{p/2}\\
 &+
 c\left( \frac{\sum_{i=n_k+1}^{n_{k+1}} \big[\big(\Sbep[Y_i]\big)^++\big(\Sbep[-Y_i]\big)^+\big]}{a_{n_k} }\right)^p\\
\le    c  \sum_{i=n_k+1}^{n_{k+1}} & \frac{ \Sbep\big[(|X_i|\wedge a_i)^p]}{a_i^p}+
 c\left( \frac{n_{k+1}-n_k}{n_k\log\log n_k}\right)^{p/2}+c
  \left( \frac{n_{k+1}-n_k}{n_k}\right)^{p}.
 \end{align*}
By Lemma \ref{lem2} (ii), it follows that
\begin{align*}
\sum_{k=1}^{\infty}\Capc\left(II_k\ge  \delta a_{n_k}\right)\le c  \sum_{i=1}^{\infty} \frac{ \Sbep\big[(|X_1|\wedge a_i)^p]}{a_i^p}+
 c\sum_{k=1}^{\infty} \left( \frac{1}{k^{\alpha}}\right)^{p/2}<\infty, \forall \delta>0,
 \end{align*}
 whenever we choose $p>2$ such that $\alpha p/2>1$.
 It follows that
\begin{equation}\label{eqLIL1.16} \outCapc\left(\Big\{\frac{II_k}{a_{n_k}}>\delta\Big\}\; i.o.\right)=0,\;\; \forall \delta.
\end{equation}

Finally, we consider the first term $I_k$. Let $y=2b_{n_k}$ and $x=(1+\epsilon)^2 a_{n_k}$. Then $|Y_i-\Sbep[Y_i]|\le y$ and
$xy\le \frac{\epsilon}{10}n_k$. By (\ref{eqth2.1}), we have
\begin{align*}
\Capc\left(I_k\ge  (1+\epsilon)^2 a_{n_k}\right)\le
\exp\left\{-\frac{(1+\epsilon)^4 a_{n_k}^2}{2\big(\epsilon n_k/10+\sum_{i=1}^{n_k}\Sbep[|Y_i-\Sbep[Y_i]|^2]\big)}\right\}.
\end{align*}
Since
$$ \left|\Sbep[X_i^2]-\Sbep[Y_i^2]\right|\le \Sbep|X_i^2-Y_i^2|=\Sbep[(X_1^2-b_i^2)^+]\to 0, \text{ as } i\to \infty$$
and $|\Sbep[Y_i]|\to 0$, $|\cSbep[Y_i]|\to 0$ as $i\to \infty$, 
we have $\sum_{i=1}^{n_k}\Sbep[Y_i^2]\le (1+\epsilon/2)n_k \Sbep X_1^2=(1+\epsilon/2)n_k$ for $k$ large enough.
It follows that
\begin{align*}
  \sum_{k=k_0}^{\infty} & \Capc\left(I_k\ge   (1+\epsilon)^2 a_{n_k}\right)\le
\sum_{k=k_0}^{\infty}\exp\left\{-(1+\epsilon)^2 \log\log n_k\right\}
   \le   \sum_{k=k_0}^{\infty}  \frac{c}{k^{(1+\epsilon)(1-\alpha)}} <\infty
\end{align*}
if $\alpha$ is chosen such that $ (1+\epsilon)(1-\alpha)>1$.  It follows that by the countably sub-additivity and the Borel-Cantelli Lemma again,
\begin{equation}\label{eqLIL1.17} \outCapc\left(\Big\{\frac{I_k}{a_{n_k}}>(1+\epsilon)^2\Big\}\; i.o.\right)=0.
\end{equation}
Combining (\ref{eqLIL1.15})-(\ref{eqLIL1.17}) yields (\ref{eqLIL1.12}).  The proof of (\ref{eqLIL1.10}) is proved.

From (\ref{eqLIL1.10}), it follows that
$$ \liminf_{n\to \infty}\frac{S_n}{a_n}=-\limsup_{n\to \infty}\frac{-S_n}{a_n}\ge  -1\;\; a.s.\outCapc.$$
The part (a) is proved.

\noindent{\bf Step 2}. We show (\ref{eqthLIL.3})$\implies$(\ref{eqthLIL.1}).  Suppose $C_{\Capc}\left[\frac{X_1^2}{\log\log|X_1|}\right]=\infty$.
Then, by (\ref{eq1.3}) and Lemma \ref{lem2} (i),
\begin{align}\label{eqthLIL1.20}
 \sum_{j=1}^{\infty} \Capc\big(|X_j|> M(1-\epsilon) a_j)
  \ge & \sum_{j=1}^{\infty}\Sbep\left[g_{\epsilon}\big(\frac{|X_j|}{M a_j}\big)\right]
= \sum_{j=1}^{\infty}\Sbep\left[g_{\epsilon}\big(\frac{|X_1|}{M a_j}\big)\right]
\\
\ge &\sum_{j=1}^{\infty} \Capc\big(|X_1|> M a_j)=\infty, \;\; \forall M>0.\nonumber
\end{align}
If $\{X_j;j\ge 1\}$ were independent under $\Capc$ or $\outCapc$, then one can use the standard argument of the Borel-Cantelli lemma to show that
$$ \outCapc\big(\big\{|X_j|> M(1-\epsilon) a_j)\big\}\;\;i.o. \big)=1. $$
Now, since  the indicator functions are not in $C_{l,Lip}$,  we introduce the smoothing method.
Let $Z_j=g_{1/2}\big(\frac{|X_j|}{M a_j}\big)$. Then $0\le Z_j\le 1$ and $\Sbep[(Z_j-\Sbep[Z_j])^2]\le 2\Sbep[Z_j]$.
 Note $\cSbep[-Z_j+\Sbep[Z_j]]=0$. Then by (\ref{eqcor2.4}),
\begin{align*}
  \cCapc\Big(\sum_{j=1}^n Z_j\le & \frac{1}{2} \sum_{j=1}^n \Sbep[Z_j]\Big)=
  \cCapc\Big(\sum_{j=1}^n (-Z_j+\Sbep[Z_j])\ge \frac{1}{2}\sum_{j=1}^n\Sbep[Z_j] \big)\Big) \\
\le & C \frac{8 \sum_{j=1}^n \Sbep[Z_j]}{( \sum_{j=1}^n \Sbep[Z_j])^2}=C \frac{8  }{  \sum_{j=1}^n \Sbep[Z_j]  }\to 0.
\end{align*}
 So
$$
 \Capc\Big(\sum_{j=1}^n Z_j>\frac{1}{2} \sum_{j=1}^n \Sbep[Z_j]\Big)\to 1\;\;\text{ as }  n\to \infty.
$$

If $\Capc$ is continuous as assumed in Theorem \ref{thLIL} (b), then $\Capc\equiv \outCapc$.
 If $\Sbep$ is sub-additive as assumed in Theorem   \ref{thLIL} (c), then
\begin{equation}\label{eqLIL1.20ad}\outCapc(|X|\ge c)\le \Capc(|X|\ge c)\le\Sbep\left[g_{\delta}(|X|/c)\right]\le \outCapc\big(|X|\ge c(1-\delta)\big),\; \forall \delta>0
\end{equation}
by (\ref{eq1.3}) and (\ref{eq3.2}).  Write
$A_n=\{\sum_{j=1}^n Z_j>\frac{1}{4} \sum_{j=1}^n \Sbep[Z_j]\}$. In either case, we have
 $$
 \outCapc\left(A_n\right)
 \ge \Capc\left(\sum_{j=1}^n Z_j>\frac{1}{2} \sum_{j=1}^n \Sbep[Z_j]\right)\to 1\;\;\text{ as }  n\to \infty.
 $$

Now,
\begin{align}\label{eqthLIL1.21}
\outCapc & \left(\limsup_{n\to\infty} \frac{|X_n|}{a_n}>M(1-\epsilon)\right)
= \outCapc\left(\Big\{\frac{|X_j|}{M a_j}> (1-\epsilon)\Big\} \;\;i.o.\right)\\
& \; \ge
\outCapc\left(\sum_{j=1}^{\infty} g_{\epsilon}\big(\frac{|X_j|}{M a_j}\big)=\infty\right)
=   \outCapc\left( A_n\;\; i.o.\right)
\ge   \lim_{n\to \infty}  \outCapc\left(A_n\right)=1,
\nonumber
\end{align}
 by the continuity of $\outCapc$. On the other hand,
$$\limsup_{n\to\infty} \frac{|X_n|}{a_n}\le \limsup_{n\to\infty}\Big(\frac{|S_n|}{a_n}+\frac{|S_{n-1}|}{a_n}\Big)
\le 2\limsup_{n\to\infty} \frac{|S_n|}{a_n}.
$$
It follows that
$$ \outCapc\left(\limsup_{n\to\infty} \frac{|S_n|}{a_n}>m\right)=1,\;\; \forall m>0. $$
Hence
$$ \outCapc\left(\limsup_{n\to\infty} \frac{|S_n|}{a_n}=+\infty\right)=\lim_{m\to \infty} \outCapc\left(\limsup_{n\to\infty} \frac{|S_n|}{a_n}>m\right)=1, $$
which  contradict (\ref{eqthLIL.3}). So, (\ref{eqthLIL.1}) holds.

\noindent{\bf Step 3}. Finally, we show $\Sbep[X_1]=\Sbep[-X_1]=0$.
   If $\Sbep$ is countably additive as assumed in (b). Then (\ref{eqthLIL.1}) implies
 $\Sbep \left[\frac{X_1^2}{\log\log|X_1|}\right]\le C_{\Capc}\left[\frac{X_1^2}{\log\log|X_1|}\right]<\infty$ by Lemma~\ref{lem2} (iii).
 And  so $\lim\limits_{c\to \infty}  \Sbep[(|X_1|-c)^+]=0$. Note $\lim\limits_{c\to \infty} \Sbep[(-c)\vee(X_1\wedge c)]= \Sbep[X_1]$.
 Write $Y_j= (-c)\vee(X_j\wedge c) $.   Then $\Sbep[Y_j]=\Sbep[Y_1]\to \Sbep [X_1]$ as $c\to +\infty$. So, for $c$ large enough, by (\ref{eqcor2.4}) we have
 \begin{align*}
 \cCapc\Big(  \frac{S_n}{n}<\Sbep[X_1]-2\epsilon\Big) \le &  \cCapc\Big(  \sum_{k=1}^n \big(-Y_j+\Sbep[Y_j]\big)>n\epsilon\Big)
 +   \Capc\Big(  \sum_{k=1}^n |X_j-Y_j|>n\epsilon\Big) \\
 \le & C \frac{\sum_{k=1}^n \Sbep\big[\big(-Y_j+\Sbep[Y_j]\big)^2\big]}{n^2\epsilon^2}
 +\frac{\sum_{k=1}^n \Sbep\big[|X_j-Y_j|\big]}{n\epsilon}\\
 \le & \frac{C c^2}{n\epsilon^2}+\frac{\Sbep[(|X_1|-c)^+]}{\epsilon}\to 0\;\;\text{ as } n\to \infty \text{ and then } c\to \infty.
 \end{align*}
By (\ref{eqLIL1.20ad}), it follows that
\begin{align*} &\lim_{n\to \infty}\outCapc\left(\frac{S_n}{n}>\Sbep[X_1]-3\epsilon\right)
\ge \lim_{n\to \infty}\Capc\left(\frac{S_n}{n}>\Sbep[X_1]-2\epsilon\right)=1.
\end{align*}
By the continuity of $\outCapc$,
$$\outCapc\left( \limsup_{n\to \infty}\frac{S_n}{n}>\Sbep[X_1]-2\epsilon\right)
\ge  \limsup_{n\to \infty} \outCapc\left(\frac{S_n}{n}>\Sbep[X_1]-3\epsilon\right)=1. $$
It follows that
\begin{equation}\label{eqthLIL1.22}
\outCapc\left(\limsup_{n\to \infty}\frac{S_n }{n}\ge\Sbep[X_1]\right)
=\lim_{\epsilon\downarrow 0}\outCapc\left( \limsup_{n\to \infty}\frac{S_n}{n}>\Sbep[X_1]-3\epsilon\right)=1,
 \end{equation}
by the continuity of $\outCapc$ again.  Write
 $$ A=\left\{\limsup_{n\to \infty}\frac{S_n }{n}\ge\Sbep[X_1]\right\},\; B=\left\{\limsup_{n\to \infty}\frac{|S_n| }{a_n}=+\infty\right\}. $$
 By (\ref{eqthLIL.3}), (\ref{eqthLIL1.22}) and the sub-additivity of $\outCapc$, it follows that
 $$ \outCapc(AB^c)=\outCapc(A\setminus AB)\ge \outCapc(A)-\outCapc(AB)\ge \outCapc(A)-\outCapc(B)>0. $$
 However, on $B^c$ we have $\limsup_{n\to \infty}\frac{|S_n|}{n}=0$. It follows that
 $ \Sbep[X_1]\le 0$.  Similarly, $ \Sbep[-X_1]\le 0$. So,
 $0\le -\Sbep[-X_1]\le \Sbep[X_1]\le 0. $
 Hence $\Sbep[X_1]=\Sbep[-X_1]=0$. (b) and (c) are now proved. $\Box$

 \bigskip

 {\bf Proof of Corollary \ref{cor3.1}.}
   Note that if a sequence of $\{x_n\}$ satisfies $x_n-x_{n-1}\to 0$, then $$C\big(\{x_n\}\big)=\big[\liminf_{n\to \infty} x_n, \limsup_{n\to \infty} x_n\big]. $$
 Now, it can be showed that
 $$ \frac{S_n}{a_n}-\frac{S_{n-1}}{a_{n-1}}=\frac{X_n}{a_n}+\frac{S_{n-1}}{a_{n-1}}\left(\frac{a_{n-1}}{a_n}-1\right)\to 0\;\; a.s. \Capc. $$
So,
 \begin{equation}\label{eqproofcor3.1.1}
\cCapc\left(C\Big\{\frac{S_n}{a_n}\Big\}=\Big[\liminf_{n\to \infty}\frac{S_n}{a_n},\;\;\limsup_{n\to \infty}\frac{S_n}{a_n}\Big]\right)=1.
\end{equation}
On the other hand,   by Theorem \ref{thLIL}  we have
 \begin{equation}\label{eqproofcor3.1.2}  \cCapc\left(\limsup_{n\to \infty} \frac{|S_n|}{a_n}\le \overline{\sigma} \right)=1.
 \end{equation}
Note the facts $\cCapc\left(A\bigcap_{i=1}^m A_i\right)\ge \cCapc(A)+\sum_{i=1}^m\big( \cCapc(A_i)-1\big)$ and $\Capc\left(A\bigcap_{i=1}^m A_i\right)\ge \Capc(A)+\sum_{i=1}^m\big( \cCapc(A_i)-1\big)$. By (\ref{eqproofcor3.1.1}) and (\ref{eqproofcor3.1.2}), for (I) it is sufficient to show that
$$
\Capc\left(\liminf_{n\to \infty}\frac{S_n}{a_n}\le -\overline{\sigma}  \text{ and } \limsup_{n\to \infty}\frac{S_n}{a_n}\ge \overline{\sigma} \right)=1
$$
and, for (II) and (III) it is sufficient to show that
$$
\cCapc\left(\liminf_{n\to \infty}\frac{S_n}{a_n}\le-\underline{\sigma}\right)=1 \text{ and } \cCapc\left(\limsup_{n\to \infty}\frac{S_n}{a_n}\ge \underline{\sigma}\right)=1.
$$
By the continuity of $\Capc$ (and $\cCapc$), it is sufficient to show that
 \begin{equation}\label{eqproofcor3.1.3}
\Capc\left(\liminf_{n\to \infty}\frac{S_n}{a_n}\le-\overline{\sigma}+\epsilon \text{ and } \limsup_{n\to \infty}\frac{S_n}{a_n}\ge \overline{\sigma}-\epsilon\right)=1,\;\; \forall \epsilon>0
\end{equation}
and
 \begin{equation}\label{eqproofcor3.1.4}
\cCapc\left(\liminf_{n\to \infty}\frac{S_n}{a_n}\le -\underline{\sigma}+\epsilon\right)=1 \text{ and } \cCapc\left(\limsup_{n\to \infty}\frac{S_n}{a_n}\ge \underline{\sigma}-\epsilon\right)=1,\;\; \forall \epsilon>0.
\end{equation}

Let $n_k=k^k$. Then $n_{k-1}/n_k\to 0$ and $a_{n_{k-1}}/a_{n_k}\to 0$.  Note (\ref{eqproofcor3.1.2}) and
 \begin{equation}\label{eqproofcor3.1.5} \frac{S_{n_k}}{a_{n_k}}=\frac{S_{n_k}-S_{n_{k-1}}}{\sqrt{2(n_k-n_{k-1})\log\log n_k}}\sqrt{1-\frac{n_{k-1}}{n_k}}+
 \frac{S_{n_{k-1}}}{a_{n_{k-1}}}\frac{a_{n_{k-1}}}{a_{n_k}}.
 \end{equation}
  So, for (\ref{eqproofcor3.1.4}) it is sufficient to show that for any $b$ with $|b|< \underline{\sigma}$,
\begin{equation}\label{eqproofcor3.1.6}   \cCapc\left( \liminf_{k\to \infty}\left|\frac{S_{n_k}-S_{n_{k-1}}}{\sqrt{2(n_k-n_{k-1})\log\log n_k}}-b\right|<\epsilon\right)=1, \;\; \forall \epsilon>0.
\end{equation}
By the independence of the sequence $\{S_{n_k}-S_{n_{k-1}}; k\ge 2\}$ and the smoothing argument as showing (\ref{eqthLIL1.21})
from (\ref{eqthLIL1.20}), it is sufficient to prove
\begin{equation}\label{eqproofcor3.1.7}
\sum_{k=1}^{\infty} \cCapc\left(\left|\frac{S_{n_k}-S_{n_{k-1}}}{\sqrt{2(n_k-n_{k-1})\log\log n_k}}-b\right|<\epsilon\right)=\infty, \;\; \forall \epsilon>0.
\end{equation}
Applying Theorem \ref{thLower} (a) with $y_n= \sqrt{2\log\log n_k}$ yields
\begin{align*}
& \cCapc\left(\left|\frac{S_{n_k}-S_{n_{k-1}}}{\sqrt{2(n_k-n_{k-1})\log\log n_k}}-b\right|<\epsilon\right)\\
\ge &
\exp\left\{-\Big(\big(b/\underline{\sigma}\big)^2+\delta\Big)\log\log n_k\right\}\ge c k^{-\big((b/\underline{\sigma})^2+\delta\big)},
\end{align*}
if $k$ is large enough, where $(b/\underline{\sigma})^2+\delta<1$.   (\ref{eqproofcor3.1.7}) follows
 and the proof of (\ref{eqproofcor3.1.6}) is completed.

 For (\ref{eqproofcor3.1.3}), by applying Theorem \ref{thLower} (b) instead of Theorem \ref{thLower} (a) we have
 \begin{equation}\label{eqproofcor3.1.8}
\sum_{k=1}^{\infty} \Capc\left(\left|\frac{S_{n_k}-S_{n_{k-1}}}{\sqrt{2(n_k-n_{k-1})\log\log n_k}}-b\right|<\epsilon\right)=\infty, \;\; \forall \; |b|< \overline{\sigma}, \epsilon>0.
\end{equation}
 Let $f_{\epsilon}(x)$ be a continuous function such that $I\{x<\epsilon/2\}\le f_{\epsilon}(x)\le I\{x<\epsilon\}$.  Denote
 $$\eta_k =\frac{S_{n_k}-S_{n_{k-1}}}{\sqrt{2(n_k-n_{k-1})\log\log n_k}}. $$
 With the similar argument as showing (\ref{eqthLIL1.21})
from (\ref{eqthLIL1.20}), it can be verified that
$$ \Capc\left(\sum_{j=1}^{\infty} f_{\epsilon}(|\eta_j-b|)=\infty\right)=1, \;\; \forall \; |b|< \overline{\sigma}, \epsilon>0. $$
Now, for $b_1$ and $b_2$ with $|b_1|, |b_2|<\overline{\sigma}$ we have
\begin{align}
&\Capc\left(\sum_{j=1}^{\infty} f_{\epsilon}(|\eta_j-b_1|)=\infty \text{ and } \sum_{j=1}^{\infty} f_{\epsilon}(|\eta_j-b_2|)=\infty\right) \nonumber\\
=& \lim_{l_1\to\infty}\lim_{m_1\to \infty}  \Capc\left(\sum_{j=1}^{m_1} f_{\epsilon}(|\eta_j-b_1|)\ge l_1 \text{ and } \sum_{j=1}^{\infty} f_{\epsilon}(|\eta_j-b_2|)=\infty\right)\nonumber\\
=& \lim_{l_1\to\infty}\lim_{m_1\to \infty}   \Capc\left(\sum_{j=1}^{m_1} f_{\epsilon}(|\eta_j-b_1|)\ge l_1 \text{ and } \sum_{j=m_1+1}^{\infty} f_{\epsilon}(|\eta_j-b_2|)=\infty\right)\nonumber\\
=& \lim_{l_1\to\infty}\lim_{m_1\to \infty}  \lim_{l_2\to\infty}\lim_{m_2\to \infty}  \Capc\left(\sum_{j=1}^{m_1} f_{\epsilon}(|\eta_j-b_1|)\ge l_1 \text{ and } \sum_{j=m_1+1}^{m_2} f_{\epsilon}(|\eta_j-b_2|)\ge l_2\right)\nonumber\\
\ge & \lim_{l_1\to\infty}\lim_{m_1\to \infty}  \lim_{l_2\to\infty}\lim_{m_2\to \infty}  \Capc\left(\sum_{j=1}^{m_1} f_{\epsilon}(|\eta_j-b_1|)\ge 2 l_1\right)\nonumber\\
& \qquad\qquad \qquad\qquad\qquad \cdot\Capc\left( \sum_{j=m_1+1}^{m_2} f_{\epsilon}(|\eta_j-b_2|)\ge 2 l_2\right)\nonumber\\
=&\Capc\left(\sum_{j=1}^{\infty} f_{\epsilon}(|\eta_j-b_1|)=\infty\right)\Capc\left( \sum_{j=1}^{\infty} f_{\epsilon}(|\eta_j-b_2|)=\infty\right)=1.
\end{align}
The inequality above is due to the fact that
$$ \Capc(X\ge l_1, Y\ge l_2)\ge \Sbep[g(X/l_1)g(Y/l_2)]=\Sbep[g(X/l_1)]\Sbep[g(Y/l_2)]\ge \Capc (X\ge 2 l_1)\Capc(Y\ge 2 l_2) $$
if $Y$ is independent to  $X$ under $\Sbep$, where $g(x)$ is a continuous function such that $I\{x\ge 1\}\ge g(x)\ge I\{x\ge 2 \}$.  By noting that $\sum_{j=1}^{\infty} f_{\epsilon}(|\eta_j-b |)=\infty$ implies
$\liminf\limits_{k\to \infty}|\eta_k-b|<\epsilon$, we conclude that
$$ \Capc\left(\liminf_{k\to \infty}|\eta_k-b_1|<\epsilon \text{ and } \liminf_{k\to \infty} |\eta_k-b_2|<\epsilon\right)=1\;\; \forall  |b_1|, |b_2|<\overline{\sigma}, \epsilon>0. $$
Combining the above equality with (\ref{eqproofcor3.1.2}) and (\ref{eqproofcor3.1.5}) yields
$$ \Capc\left(\liminf_{n\to \infty}\left|\frac{S_n}{a_n}-b_1\right|<\epsilon \text{ and }
\liminf_{n\to \infty}\left|\frac{S_n}{a_n}-b_2\right|<\epsilon \right)=1\;\; \forall  |b_1|, |b_2|<\overline{\sigma}, \epsilon>0. $$
Now, (\ref{eqproofcor3.1.3}) follows by letting $b_1=-(1-\epsilon)\overline{\sigma}$ and $b_2=(1-\epsilon)\overline{\sigma}$.
 The proof of the corollary is completed.  $\Box$

\Acknowledgements{Special thanks go to the anonymous referees, the associate editor and
the editors for their constructive comments, which led to a much improved
version of this paper. This work was supported by grants from the NSF of China (No. 11225104), the 973 Program
(No. 2015CB352302) and the Fundamental Research Funds for the Central Universities.}


\end{document}